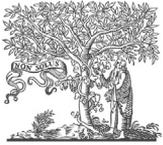

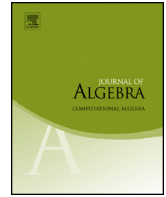

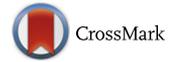

# The ideals of the slice Burnside $p$-biset functor ☆

## Ibrahima Tounkara

*Laboratoire d'Algèbre, de Cryptologie, de Géométrie Algébrique et Applications (LACGAA), Département de Mathématiques et Informatique, Faculté des Sciences et Techniques, Université Cheikh Anta Diop, BP 5005 Dakar, Senegal*

| ARTICLE INFO | ABSTRACT |
|---|---|
|  | Let $G$ be a finite group and $\mathbb{K}$ be a field of characteristic zero. Our purpose is to investigate the ideals of the slice Burnside functor $\mathbb{K}\Xi$. It turns out that they are the subfunctors $F$ of $\mathbb{K}\Xi$ such that for any finite group $G$, the evaluation $F(G)$ is an ideal of the algebra $\mathbb{K}\Xi(G)$. This allows for a determination of the full lattice of ideals of the slice Burnside $p$-biset functor $\mathbb{K}\Xi_p$.<br><br>© 2017 Published by Elsevier Inc. |

## 1. Introduction

The biset category $\mathcal{C}$ of finite groups has all finite groups as objects, the group of morphisms from a finite group $G$ to a finite group $H$ is the double Burnside group $B(H, G)$, i.e. the Grothendieck group of $(H, G)$-bisets. In particular the endomorphism ring of a finite group $G$ is the double Burnside ring $B(G, G)$.

A biset functor is an additive functor over this preadditive category, with values in abelian groups, and biset functors form an abelian category $\mathcal{F}$. More generally, one can extend the morphisms to $RB(H, G) = R \otimes_{\mathbb{Z}} B(H, G)$, where $R$ is a commutative ring, and consider the $R$-linear functors with values in the category of $R$-modules. Thus, one obtain an $R$-linear abelian category $\mathcal{F}_R$.

A fundamental example of biset functor is the Burnside functor: it can be viewed as the representable functor $RB(1, -)$, or the Yoneda functor corresponding to the trivial group.

In particular it is a projective object of the category $\mathcal{F}_R$, which allows by the Yoneda–Dress construction to build enough projective objects on the category $\mathcal{F}_R$.

Moreover, the Burnside functor has a multiplicative structure which endows it with a structure of Green biset functor, and the modules over this Green functor are the biset functors.

We have a good parametrization of the isomorphism classes of simple $RB$-modules by isomorphism classes of pairs $(H, V)$, where $H$ is a finite group and $V$ is a simple $R\mathrm{Out}(H)$-module. To each such $(H, V)$ corresponds the isomorphism class of $S_{H,V}$, where $S_{H,V}(G)$ is the quotient of $L_{H,V}(G) = RB(G, H) \otimes_{RB(H,H)} V$ by

$$J_{H,V}(G) = \{\sum_i \phi_i \otimes v_i \in L_{H,V}(G) \mid \forall \psi \in \mathbb{K}B(H, G), \ \sum_i (\psi\phi_i).v_i = 0\}.$$

An important property of $H$ is that it is a minimal group for $S_{H,V}$, and all the minimal groups for a simple biset functor are isomorphic.

Note that in general, the explicit computation of the evaluation $S_{H,V}(G)$ of a simple functor or more generally of a simple module over a Green biset functor is not easy (cf. [9] and [5]). The study of the functor $\mathbb{K}B$ (cf. [1]) where $\mathbb{K}$ is a field of characteristic zero, has allowed for an explicit description of some simple biset functors by the introduction of a new class of finite groups, the **B**-groups.

In this paper, we consider the slice Burnside ring introduced in ([3]). It is an analogue of the classical Burnside ring constructed from the morphisms of $G$-sets instead the $G$-sets themselves, and it shares most of its properties. In particular, as already shown by Serge Bouc (see [3] for more complete description), the slice Burnside ring is a commutative ring, which is free of finite rank as a $\mathbb{Z}$-module, and it becomes a split semisimple $\mathbb{Q}$-algebra, after tensoring with $\mathbb{Q}$. The correspondence which assigns to each finite group its slice Burnside ring has a natural biset functor structure, for which it becomes a Green biset functor.

We restrict our attention to the subcategory $\mathcal{C}_p$ of $\mathcal{C}$ whose objects are all the finite $p$-groups (for a given prime number $p$). To achieve the description of the ideals of the slice Burnside $p$-biset functor over a field of characteristic zero, we introduce the concept of **T**-slices instead of **B**-groups. In studying the ideals of the slice Burnside $p$-biset functor, we find another counterexample to Serge Bouc's conjecture saying that the minimal groups for a simple module over a Green biset functor form a single isomorphism class.

The study in this direction has it origin in Nadia Romero's paper ([8]), in which she found the first counterexample to this conjecture for the monomial Burnside ring over a cyclic group of order four.

The plan of the paper is as follows. In Section 2 we recall the basic definitions and preliminary properties of the classical Burnside functor. Section 3 recalls some properties of the slice Burnside ring and the slice Burnside functor. Section 4 is devoted to the action of biset operations on idempotents. In that section we exploit the case of the deflation and find a constant $m_{G,S,N}$. Section 5 gives some properties of the constant $m_{G,S,N}$. Section 6, exploits $m_{G,S,N}$ in the case of $p$-groups. Section 7 gives a characterization of the ideals of the slice Burnside functor. In Section 8 we further specialize the biset category to be the $p$-biset category and we consider the slice Burnside $p$-biset functor. In the final Section, we give a counterexample to Bouc's conjecture.

## 2. Definitions and preliminary properties

The aim of this section is to fix some notation and recall some properties of biset functors. Throughout this paper, we fix a field $\mathbb{K}$ of characteristic zero and we identify the prime subfield of $\mathbb{K}$ with $\mathbb{Q}$. We assume that the reader is familiar with the concept of biset.

**Notation 2.1.** If $Y$ is an $(H, G)$-biset and $X$ is a $(K, H)$-biset then the *tensor product* $X \times_H Y$ is defined as the set of $H$-orbits of $X \times Y$ under the $H$-action defined by $(x, y).h := (xh, h^{-1}.y)$ for $h \in H$ and $(x, y) \in X \times Y$. The $H$-orbit of $(x, y)$ is denoted by $(x, {}_H y)$ and the set $X \times_H Y$ of elements $(x, {}_H y)$ is a $(K, G)$-biset under

$$k(x, {}_H y)g := (kx, {}_H yg).$$

This operation defines a bilinear map

$$- \times_H - : B(K, H) \times B(H, G) \to B(K, G), \ ([X], [Y]) \mapsto [X \times_H Y],$$

where $B(H, G)$ is the Burnside group of $(H, G)$-bisets. If $G$ is a finite group, we denote by $\mathrm{Id}_G$ the $(G, G)$-biset $G$.

This leads to the following definitions.

**Definition 2.2.** The *biset category* $\mathcal{C}$ of finite groups is the category defined as follows:

- The objects of $\mathcal{C}$ are finite groups.
- If $G$ and $H$ are finite groups, then $\mathrm{Hom}_{\mathcal{C}}(G, H) = B(H, G)$.
- If $G$, $H$, and $K$ are finite groups, then the composition $v \circ u$ of the morphism $u \in \mathrm{Hom}_{\mathcal{C}}(G, H)$ and the morphism $v \in \mathrm{Hom}_{\mathcal{C}}(H, K)$ is equal to $v \times_H u$.
- For any finite group $G$, the identity morphism of $G$ in $\mathcal{C}$ is equal to $\mathrm{Id}_G$.

More generally, when $\mathbb{K}$ is a commutative ring, we define similarly the category $\mathbb{K}\mathcal{C}$ by extending coefficients to $\mathbb{K}$, i.e. by setting

$$\operatorname{Hom}_{\mathbb{K}\mathcal{C}}(G,H) = \mathbb{K} \otimes_{\mathbb{Z}} B(H,G),$$

which will be simply denoted by $\mathbb{K}B(H,G)$.

**Remark 2.3.** The category $\mathcal{C}$ is a preadditive category: the sets of morphisms in $\mathcal{C}$ are abelian groups, and the composition of morphisms is bilinear. Let $G$ and $H$ be finite groups. Then, any morphism from $G$ to $H$ in $\mathcal{C}$ is a linear combination with integral coefficients of morphisms of the form $[(G \times H)/L]$, where $L$ is some subgroup of $H \times G$.

**Notation 2.4.** Let $H$ and $G$ be finite groups.

- We indicate by $H \leq G$ that $H$ is a subgroup of $G$. We write $H < G$ if $H \leq G$ and $H \neq G$.
- We set $H^x := x^{-1}Hx$, for $x \in G$ and $H \leq G$.
- If $H \leq G$, we denote by $\operatorname{Res}_H^G$ the set $G$, viewed as an $(H,G)$-biset for left and right multiplication, and by $\operatorname{Ind}_H^G$ the same set viewed as a $(G,H)$-biset.
- If $N \trianglelefteq G$, and $H = G/N$, we denote by $\operatorname{Inf}_H^G$ the set $H$, viewed as a $(G,H)$-biset for the left action of $G$, and right action of $H$ by multiplication. Also we denote by $\operatorname{Def}_H^G$ the set $H$, viewed as an $(H,G)$-biset.
- If $f : G \to H$ is a group isomorphism, we denote by $\operatorname{Iso}_G^H$ or $\operatorname{Iso}(f)$ the set $H$, viewed as an $(H,G)$-biset for left multiplication in $H$, and right action of $G$ given by multiplication by the image under $f$.
  The bisets $\operatorname{Ind}_H^G$, $\operatorname{Res}_H^G$, $\operatorname{Inf}_H^G$, $\operatorname{Def}_H^G$ and $\operatorname{Iso}_G^H$ are called the *elementary bisets*.

**Lemma 2.5.** *[1.1.3 Relations [2]] Commutation conditions:*

- *(Mackey formula) If $H$ and $K$ are subgroups of $G$, then*

$$\operatorname{Res}_H^G \circ \operatorname{Ind}_K^G \cong \sum_{x \in [H\backslash G/K]} \operatorname{Ind}_{H \cap {}^xK}^H \circ \operatorname{Iso}(\gamma_x) \circ \operatorname{Res}_{H^x \cap K}^K,$$

  *where $[H\backslash G/K]$ is a set of representatives of $(H,K)$-double cosets in $G$, and $\gamma_x$ is the group isomorphism induced by conjugation by $x$.*
- *If $H$ is a subgroup of $G$, and if $N$ is a normal subgroup of $G$, then*

$$\operatorname{Def}_{G/N}^G \circ \operatorname{Ind}_H^G \cong \operatorname{Ind}_{HN/N}^{G/N} \circ \operatorname{Iso}(\varphi) \circ \operatorname{Def}_{H/H\cap N}^H,$$

  *where $\varphi : H/H \cap N \to HN/N$ is the canonical group isomorphism.*

The following is a straightforward consequence of Goursat's Lemma on subgroups of a direct product of two groups:

**Lemma 2.6.** *[[2], Lemma 2.3.26] Any transitive $(H, G)$-biset is isomorphic to a composition*

$$\mathrm{Ind}_D^H \circ \mathrm{Inf}_{D/C}^D \circ \mathrm{Iso}_{D/C}^{B/A} \circ \mathrm{Def}_{B/A}^B \circ \mathrm{Res}_B^G,$$

*where $A \trianglelefteq B \leq G$, $C \trianglelefteq D \leq H$, and $f : B/A \xrightarrow{\simeq} D/C$ is a group isomorphism.*

**Remark 2.7.** Let $\mathcal{E}$ be the set of triples $\big((D, C), f, (B, A)\big)$ where $A \trianglelefteq B \leq G$, $C \trianglelefteq D \leq H$, and $f : B/A \xrightarrow{\simeq} D/C$ is a group isomorphism. The group $H \times G^{op}$ (where $G^{op}$ is the opposite group of $G$) acts by conjugation on $\mathcal{E}$ and the set of all elements $\mathrm{Ind}_D^H \circ \mathrm{Inf}_{D/C}^D \circ \mathrm{Iso}_{D/C}^{B/A} \circ \mathrm{Def}_{B/A}^B \circ \mathrm{Res}_B^G$ is a $\mathbb{K}$-basis of $\mathbb{K}B(H, G)$, where the triple $\big((D, C), f, (B, A)\big)$ runs over representatives of $(H \times G)$-orbits in $\mathcal{E}$.

**Definition 2.8.** A *biset functor* is an additive functor from $\mathcal{C}$ to abelian groups. A biset functor with values in $\mathbb{K}$-Vect is a $\mathbb{K}$-linear functor from $\mathbb{K}\mathcal{C}$ to the category $\mathbb{K}$-Vect of $\mathbb{K}$-vector spaces.

**Remark 2.9.** Biset functors form an abelian category $\mathcal{F}$ (where morphisms are natural transformations of functors).

Similarly, biset functors with values in $\mathbb{K}$-Vect form a $\mathbb{K}$-linear abelian category $\mathcal{F}_{\mathbb{K}}$.

**Definition 2.10.** Let $F$ be a biset functor on $\mathcal{C}$.

A *minimal group* for $F$ is a finite group $H$ such that $F(H) \neq 0$, but $F(K) = 0$ for any finite group $K$ with $|K| < |H|$. The class of minimal groups for $F$ is denoted by $\mathrm{Min}(F)$. Note that $\mathrm{Min}(F) \neq \emptyset$ if and only if $F$ is not the zero functor.

**Definition 2.11.** A biset functor $A$ is a *Green biset functor* (on $\mathcal{C}$) if it is endowed with bilinear products $A(G) \times A(H) \to A(G \times H)$, denoted by $(a, b) \mapsto a \times b$, for any finite groups $G$, $H$, and an element $\epsilon_A \in A(1)$, satisfying the following conditions:

1. (Associativity) Let $G$, $H$ and $K$ be finite groups. If

$$\alpha_{G,H,K} : G \times (H \times K) \to (G \times H) \times K$$

is the canonical group isomorphism, then for any $a \in A(G)$, $b \in A(H)$, and $c \in A(K)$

$$(a \times b) \times c = \mathrm{Iso}(\alpha_{G,H,K})\big(a \times (b \times c)\big).$$

2. (Identity element) Let $G$ be a finite group. Let $\lambda_G : 1 \times G \to G$ and $\rho_G : G \times 1 \to G$ denote the canonical group isomorphisms. Then for any $a \in A(G)$

$$a = \mathrm{Iso}(\lambda_G)(\epsilon_A \times a) = \mathrm{Iso}(\rho_G)(a \times \epsilon_A).$$

3. (Functoriality) If $\varphi : G \to G'$ and $\psi : H \to H'$ are morphisms in $\mathbb{K}\mathcal{C}$, then for any $a \in A(G)$ and $b \in A(H)$

$$A(\varphi \times \psi)(a \times b) = A(\varphi)(a) \times A(\psi)(b).$$

**Definition 2.12.** Let $A$ be a Green biset functor (resp. with values in $\mathbb{K}$) on $\mathcal{C}$. A *left A-module* $M$ is an object of $\mathcal{F}$ (resp. $\mathcal{F}_\mathbb{K}$) endowed, for any finite groups $G$ and $H$, with bilinear product maps

$$A(G) \times M(H) \to M(G \times H),$$

denoted by $(a, m) \to a \times m$, fulfilling the following conditions:

- For any $a \in A(G)$, $b \in A(H)$, and $m \in M(K)$

$$(a \times b) \times m = \mathrm{Iso}(\alpha_{G,H,K})\big(a \times (b \times m)\big).$$

- For any $m \in M(G)$

$$m = \mathrm{Iso}(\lambda_G)(\epsilon_A \times m).$$

- For any $a \in A(G)$ and $m \in M(H)$

$$M(\phi \times \psi)(a \times m) = A(\phi)(a) \times M(\psi)(m).$$

**Example 2.13.** Assigning to each finite group $G$ its Burnside ring $B(G)$ yields an example of a Green biset functor. The product is induced by the bifunctor sending a $G$-set $X$ and an $H$-set $Y$ to the $(G \times H)$-set $X \times Y$. The identity element is $1 \in B(1)$, i.e. the class of a set of cardinality one in $B(1) \simeq \mathbb{Z}$.

**Definition 2.14.** Let $A$ be a Green biset functor on $\mathcal{C}$. A *left ideal* of $A$ is an $A$-submodule of the left $A$-module $A$. In other words it is a biset subfunctor $I$ of $A$ such that

$$A(G) \times I(H) \subseteq I(G \times H)$$

for any finite groups $G$ and $H$.

One defines similarly a *right ideal* of $A$.

A *two sided ideal* of $A$ is a left ideal which is also a right ideal.

A Green biset functor $A$ is called *simple* if its only two sided ideals are 0 and $A$.

**Proposition 2.15** *(See [7]). Let $A$ be a Green biset functor on $\mathcal{C}$ (resp. with values in $\mathbb{K}$-Vect) and let $G$ be a finite group. Then the product*

$$. : A(G) \times A(G) \to A(G)$$
$$(u, v) \mapsto u.v = \mathrm{Iso}^G_{\triangle(G)} \mathrm{Res}^{G \times G}_{\triangle(G)}(u \times v)$$

*where △(G) denoted the diagonal subgroup of G, endows A(G) with a ring structure (resp. a structure of $\mathbb{K}$-algebra).*

*Moreover, if H is a finite group, then*

$$a \times b = (\mathrm{Inf}^{G \times H}_G a).(\mathrm{Inf}^{G \times H}_H b)$$

*for any $a \in A(G)$ and $b \in A(H)$.*

**Proposition 2.16.** *Let A be a Green biset functor. The following assertions are equivalent.*

1. *F is a subfunctor of A and the evaluation at any finite group G of F is an ideal of $A(G)$.*
2. *F is an ideal of A as Green biset functor.*

**Proof.** Let $G$ be a finite group.

Assume that Assertion (2) holds for $F$ i.e. $F$ is an ideal of $A$. Then for any $u \in A(G)$ and any $v \in F(G)$, the product $u \times v$ is in $F(G \times G)$.

Now the product $(u, v) \mapsto u.v$ in $A(G)$ can be recovered from the product $(u, v) \mapsto u \times v$ as

$$u.v = \mathrm{Iso}^G_{\triangle(G)} \mathrm{Res}^{G \times G}_{\triangle(G)}(u \times v).$$

Since $u \times v$ is in $F(G \times G)$, we have $\mathrm{Res}^{G \times G}_{\triangle(G)}(u \times v) \in F\big(\triangle(G)\big)$ and $\mathrm{Iso}^G_{\triangle(G)} \mathrm{Res}^{G \times G}_{\triangle(G)}(u \times v) \in F(G)$. Thus the product $u.v$ is in $F(G)$ and $F(G)$ is an ideal of $A(G)$.

Conversely, assume that $F$ is a subfunctor of $A$ such that $F(G)$ is an ideal of $A(G)$ for any finite group $G$.

Let $H$ be finite group.

Let $a$ be an element of $A(G)$ and let $b$ be an element of $F(H)$. Then $\mathrm{Inf}^{G \times H}_{G \times H/1 \times H}(a)$ is in $A(G \times H)$ and $\mathrm{Inf}^{G \times H}_{G \times H/G \times 1}(b)$ is in $F(G \times H)$.

Since $F(G \times H)$ is an ideal of $A(G \times H)$, we have that the product $a \times b = (\mathrm{Inf}^{G \times H}_G a).(\mathrm{Inf}^{G \times H}_H b)$ lies in $F(G \times H)$. Hence, the subfunctor F is an ideal of $A$. This completes the proof. $\square$

## 3. Review of slice Burnside ring and slice Burnside functor

We first recall the definition and basic properties of the slice Burnside ring introduced in [3], to which we refer the reader for all statements without proof.

**Definition 3.1.** A *slice* $(T, S)$ is a pair of finite groups with $S \leq T$. When $G$ is a finite group, a slice of $G$ is a slice $(T, S)$ with $T \leq G$.

**Notation 3.2.** Let $H$ and $G$ be finite groups.

- The set of slices of $G$ is denoted by $\Pi(G)$.
- When $(T, S)$ in an element of $\Pi(G)$, denote by $G/S \to G/T$ the projection morphism.
- We say that two slices $(T, S)$ and $(V, U)$ of $G$ are conjugate, if there exists an element $g$ of $G$ such that $T = V^g$ and $S = U^g$.
  We note $N_G(T, S) = N_G(T) \cap N_G(S) = \{g \in G \mid T^g = T \text{ and } S^g = S\}$.
  We write $(T, S) =_G (V, U)$ if the slices $(T, S)$ and $(V, U)$ of $G$ are conjugate, and $(T, S) \neq_G (V, U)$ otherwise.
- We say that $(V, U)$ is a quotient of $(T, S)$ and we denote $(T, S) \twoheadrightarrow (V, U)$, if there exists a surjective group homomorphism $\varphi : T \to V$ such that $\varphi(S) = U$. If $\phi$ is an isomorphism, we say that $(V, U)$ and $(T, S)$ are isomorphic.

**Definition 3.3.** Let $G$-**set** be the category of finite $G$-sets. Then the category $G$-**Mor** of morphisms of $G$-**set** has as objects the morphisms of $G$-**set**, and a morphism from $f : A \to B$ to $g : A' \to B'$ is a pair of morphisms of $G$-sets $h : A \to A'$ and $k : B \to B'$ making the following diagram commute

$$
\begin{array}{ccc}
A & \xrightarrow{\ f\ } & B \\
\downarrow{\scriptstyle h} & & \downarrow{\scriptstyle k} \\
A' & \xrightarrow{\ g\ } & B'
\end{array}
$$

Note that the category $G$-**Mor** admits products (induced by the direct product of $G$-sets) and coproducts (induced by the disjoint union of $G$-sets).

**Definition 3.4.** Let $G$ be a finite group. The *slice Burnside group* $\Xi(G)$ of $G$ is the Grothendieck group of the category of morphisms of finite $G$-sets defined as the quotient of the free abelian group on the set of isomorphism classes $[\ X \xrightarrow{\ f\ } Y\ ]$ of morphisms of finite $G$-sets, by the subgroup generated by elements of the form

$$
[\ X_1 \sqcup X_2 \xrightarrow{\ f_1 \sqcup f_2\ } Y\ ] - [\ X_1 \xrightarrow{\ f_1\ } f(X_1)\ ] - [\ X_2 \xrightarrow{\ f_2\ } f(X_2)\ ],
$$

whenever $X \xrightarrow{\ f\ } Y$ is a morphism of finite $G$-sets with a decomposition $X = X_1 \sqcup X_2$ as a disjoint union of $G$-sets, where $f_1 = f_{|X_1}$ and $f_2 = f_{|X_2}$.

The product of morphisms induces a commutative unital ring structure on $\Xi(G)$. The identity element for multiplication is the image of the class $[\bullet \to \bullet]$, where $\bullet$ denotes a $G$-set of cardinality 1. For a morphism of $G$-sets $f : X \to Y$, let $\pi(f)$ denote the image in $\Xi(G)$ of the isomorphism class of $f$.

**Theorem 3.5.** *[[3], Theorem 3.9] Let $G$ and $H$ be finite groups, and let $U$ be a finite $(H, G)$-biset. Then the correspondence*

$$(X \xrightarrow{\ f\ } Y) \ \mapsto \ (U \times_G X \xrightarrow{\ U \times_G f\ } U \times_G Y)$$

*is a functor from $G$-**Mor** to $H$-**Mor** and induces a group homomorphism $\Xi(U) : \Xi(G) \to \Xi(H)$.*

**Proof.** Since $f : X \to Y$ is a morphism of $G$-sets, the map given by

$$U \times_G f : U \times_G X \to U \times_G Y, \quad (u, {}_G x) \mapsto (u, {}_G f(x))$$

is a well defined morphism of $H$-sets. Let $f' : X' \to Y'$ be a morphism of $G$-sets and let $(\alpha, \beta)$ be a morphism from $f$ to $f'$ i.e. a pair of morphisms of $G$-sets such that $\beta \circ f = f' \circ \alpha$. Then $U \times_G f$, $U \times_G f'$, $U \times_G \alpha$ and $U \times_G \beta$ are morphisms of $H$-sets and they satisfy $(U \times_G \beta) \circ (U \times_G f) = (U \times_G f') \circ (U \times_G \alpha)$. Indeed, for any $x \in X$ and $u \in U$, we have

$$(U \times_G \beta) \circ (U \times_G f)(u, {}_G x) = \Big(u, {}_G \beta\big(f(x)\big)\Big) = \Big(u, {}_G f'\big(\alpha(x)\big)\Big)$$
$$= (U \times_G f') \circ (U \times_G \alpha)(u, {}_G x).$$

Thus $(U \times_G \alpha, U \times_G \beta)$ is a morphism from $U \times_G f$ to $U \times_G f'$ and it follows that the correspondence

$$(X \xrightarrow{\ f\ } Y) \ \mapsto \ (U \times_G X \xrightarrow{\ U \times_G f\ } U \times_G Y)$$

define a functor from $G$-**Mor** to $H$-**Mor**.

To prove to second part of the Theorem the only thing to check is that the defining relations of $\Xi(G)$ are mapped to relations in $\Xi(H)$. But if $X_1 \sqcup X_2 \xrightarrow{\ f_1 \sqcup f_2\ } Y$ is a morphism of finite $G$-sets, then $U \times_G (X_1 \sqcup X_2)$ is isomorphic to $(U \times_G X_1) \sqcup (U \times_G X_2)$.

Moreover, the image of the map $U \times_G f_1$ is equal to $U \times_G f_1(X_1)$. It follows that the relation $[\ X_1 \sqcup X_2 \xrightarrow{\ f_1 \sqcup f_2\ } Y] - [X_1 \xrightarrow{\ f_1\ } f_1(X_1)] - [X_2 \xrightarrow{\ f_2\ } f_2(X_2)]$ in $\Xi(G)$ is mapped to the relation

$$[(U \times_G f_1) \sqcup (U \times_G f_2)] - [U \times_G f_1] - [U \times_G f_2]. \qquad \square$$

**Proposition 3.6.** *The induction functor is left adjoint to the restriction functor. The adjunction means, there is a natural bijection*

$$\mathrm{Hom}_{G-\mathbf{Mor}}\big(\ \mathrm{Ind}_H^G X \xrightarrow{\ \mathrm{Ind}_H^G f\ } \mathrm{Ind}_H^G Y \ , \ A \xrightarrow{\ \alpha\ } B \ \big)$$

$$\cong \mathrm{Hom}_{H-\mathbf{Mor}}\big(\ X \xrightarrow{\ f\ } Y \ , \ \mathrm{Res}_H^G A \xrightarrow{\ \mathrm{Res}_H^G \alpha\ } \mathrm{Res}_H^G B \ \big).$$

**Proof.** There is an adjunction $\mathrm{Hom}_{G-\mathbf{set}}(\mathrm{Ind}_H^G X, A) \cong \mathrm{Hom}_{H-\mathbf{set}}(X, \mathrm{Res}_H^G A)$ which assigns to a morphism of $H$-sets $k : X \to A$ the morphism of $G$-sets $G \times_H X \to A$, $(g, {}_H x) \mapsto gk(x)$. It induces obvious bijection between $\mathrm{Hom}_{G-\mathbf{Mor}}\big(\mathrm{Ind}_H^G X \xrightarrow{\mathrm{Ind}_H^G f} \mathrm{Ind}_H^G Y$, $A \xrightarrow{\alpha} B\big)$ and $\mathrm{Hom}_{H-\mathbf{Mor}}\big(X \xrightarrow{f} Y$, $\mathrm{Res}_H^G A \xrightarrow{\mathrm{Res}_H^G \alpha} \mathrm{Res}_H^G B\big)$. $\quad\square$

Similarly,

**Proposition 3.7.** *The deflation functor is left adjoint to the inflation functor in the category of morphisms of $G$-sets.*

**Proposition 3.8.** *The assignment $G \mapsto \Xi(G)$ is a biset functor.*

**Proof.** If $G$ and $H$ are finite groups, and if $U$ is a finite $(H, G)$-biset, then the functor

$$I_U : (X \xrightarrow{f} Y) \mapsto (U \times_G X \xrightarrow{U \times_G f} U \times_G Y)$$

induces an abelian group homomorphism $\Xi(U) : \Xi(G) \to \Xi(H)$, by Theorem 3.5. If $U'$ is an $(H, G)$-biset isomorphic to $U$, then the functors $I_U$ and $I_{U'}$ are clearly isomorphic. Hence $\Xi(U) = \Xi(U')$. And if $U$ is the disjoint union of two $(H, G)$-bisets $U_1$ and $U_2$, then the functor $I_U$ is isomorphic to the disjoint union of the functors $I_{U_1}$ and $I_{U_2}$. It follows that $\Xi(U) = \Xi(U_1) + \Xi(U_2)$. If $K$ is another finite group, and if $V$ is a finite $(K, H)$-biset, then there is a natural isomorphism of $K$-sets $V \times_H (U \times_G X) \to (V \times_H U) \times_G X$. This induces an isomorphism of functors $I_V \circ I_U \cong I_{V \times_H U}$ and so $\Xi(V) \circ \Xi(U) = \Xi(V \times_H U)$. Finally, if $U$ is the identity biset $\mathrm{Id}_G$, then the functor $I_{\mathrm{Id}}$ is isomorphic to the identity functor. Thus $\Xi(\mathrm{Id}_G) = \mathrm{Id}_{\Xi(G)}$. $\quad\square$

**Proposition 3.9.** *[[3], Theorem 3.9] Assigning to each finite group $G$ the $\mathbb{Z}$-module $\Xi(G)$ defines a Green biset functor.*

**Proof.** We have already seen that the correspondence $G \mapsto \Xi(G)$ is a biset functor. Let $f : X \to Y$ be a morphism of $G$-sets and let $g : Z \to T$ be a morphism of $H$-sets then $f \times g : X \times Z \to Y \times T$ is a morphism of $G \times H$-sets. It induces the product

$$\times : \Xi(G) \times \Xi(H) \to \Xi(G \times H)$$
$$(a, b) \mapsto a \times b.$$

Moreover, the morphism $\bullet \to \bullet$ of 1-sets is obviously an identity element for this product, up to identification $G \times 1 = G$.

Finally if $G$, $H$, $G'$, $H'$ are finite groups, if $U$ is a finite $(H, G)$-biset, if $U'$ is a finite $(H', G')$-biset, it is clear that the morphisms $(U \times U') \times_{G \times G'} (f \times f')$ and $(U \times_G f) \times (U' \times_{G'} f')$ are isomorphic morphisms of $(H \times H')$-sets.

Thus we have shown that $\Xi(G)$ satisfies the three conditions of Definition 2.11, with $\epsilon_\Xi$ being the identity element of the product in $\Xi(1)$.

So, the correspondence $G \mapsto \Xi(G)$ is a Green biset functor.  □

If $S \leq T$ are subgroups of $G$, set

$$\langle T, S \rangle_G = \pi(G/S \to G/T).$$

**Lemma 3.10.** *[[3], Lemma 3.4] Let $f : X \to Y$ be a morphism of $G$-sets. Then in the group $\Xi(G)$,*

$$\pi(f) = \sum_{x \in [G \backslash X]} \langle G_{f(x)}, G_x \rangle_G,$$

*where $G_\bullet$ denotes the stabilizer of $\bullet$.*

Thus the group $\Xi(G)$ is, as an additive group, free with basis $\{\langle T, S \rangle_G \mid (T, S) \in [\Pi(G)]\}$, where $[\Pi(G)]$ is a set of representatives of conjugacy classes of slices of $G$.

**Remark 3.11.** In this basis, the multiplication of $\Xi(G)$ can be computed as follows

$$\langle T, S \rangle_G . \langle Y, X \rangle_G = \sum_{g \in [S \backslash G / X]} \langle T \cap {}^g Y \ , S \cap {}^g X \rangle_G$$

for any slices $(T, S)$ and $(Y, X)$ of $G$.

For any slice $(T, S)$ of $G$ and any morphism of $G$-sets $f : X \to Y$, we write

$$\phi^G_{T,S}( \ X \xrightarrow{\ f\ } Y \ ) := |\mathrm{Hom}_{G-\mathbf{Mor}}(G/S \to G/T \ , \ X \xrightarrow{\ f\ } Y \ )|.$$

Then $\phi^G_{T,S}$ induces a ring homomorphism $\Xi(G) \to \mathbb{Z}$, still denoted by $\phi^G_{T,S}$. Conjugate slices give the same homomorphism.

Let

$$\Phi = \prod_{(T,S) \in [\Pi(G)]} \phi^G_{T,S} : \Xi(G) \to \prod_{(T,S) \in [\Pi(G)]} \mathbb{Z}.$$

The ring homomorphism $\Phi$ is a monomorphism between free $\mathbb{Z}$-modules having the same rank. Hence tensoring with $\mathbb{Q}$ gives a $\mathbb{Q}$-algebra isomorphism

$$\mathbb{Q} \otimes_{\mathbb{Z}} \Phi : \mathbb{Q} \otimes_{\mathbb{Z}} \Xi(G) \to \prod_{(T,S) \in [\Pi(G)]} \mathbb{Q}.$$

The commutative $\mathbb{Q}$-algebra $\mathbb{Q}\Xi(G)$ $(= \mathbb{Q} \otimes_{\mathbb{Z}} \Xi(G))$ is split semisimple.

**Notation 3.12.** For a slice $(T, S)$ of the group $G$, set

$$\xi_{T,S}^G = \frac{1}{|N_G(T,S)|} \sum_{U \leq S \leq V \leq T} |U| \mu(U,S) \mu(V,T) \langle V, U \rangle_G$$

where $\mu$ is the Möbius function of the poset of subgroups of $G$ and $N_G(T,S) = N_G(T) \cap N_G(S)$.

Note that $\xi_{T,S}^G = \xi_{T^g,S^g}^G$ for any $g \in G$.

**Proposition 3.13.** *[[3], Theorem 5.2] Let $G$ be a finite group. Then the elements $\xi_{T,S}^G$, for $(T, S) \in [\Pi(G)]$ are the primitive idempotents of $\mathbb{Q}\Xi(G)$.*

**Proof.** An easy computation shows that if $(Y, X)$ is a slice of $G$, then $\mathbb{Q}\phi_{Y,X}^G(\xi_{T,S}^G)$ is equal to 0 if $(Y, X)$ and $(T, S)$ are not $G$-conjugate, and to 1 otherwise.  □

**Remark 3.14.** After extending the scalars from $\mathbb{Q}$ to $\mathbb{K}$, we get that the commutative algebra

$$\mathbb{K}\Xi(G) = \mathbb{K} \otimes_{\mathbb{Q}} \mathbb{Q}\Xi(G)$$

is also split semisimple. We abuse notation and denote by $\xi_{T,S}^G$ the element $1 \otimes \xi_{T,S}^G$ of $\mathbb{K}\Xi(G)$.

**Proposition 3.15.** *Let $G$ be a finite group. Let $(T, S)$ be a slice of $G$. Then for any element $v$ of $\mathbb{K}\Xi(G)$*

$$v.\xi_{T,S}^G = \phi_{T,S}^G(v)\xi_{T,S}^G.$$

*Conversely, if $u \in \mathbb{K}\Xi(G)$ is such that $v.u$ is a scalar multiple of $u$ for any $v \in \mathbb{K}\Xi(G)$, then there exists $(T, S) \in \Pi(G)$ such that $u \in \mathbb{K}\xi_{T,S}^G$ and $v.u = \phi_{T,S}^G(v)u$.*

**Proof.** The set of elements $\xi_{T,S}^G$, for $(T, S) \in [\Pi(G)]$ is a basis of $\mathbb{Q}\Xi(G)$. Thus for any $v \in \mathbb{K}\Xi(G)$, there are scalars $r_{T,S}$, for $(T, S) \in \Pi(G)$, such that $v = \sum_{(T,S) \in [\Pi(G)]} r_{T,S}\xi_{T,S}^G$.

Moreover, for a slice $(V, U) \in [\Pi(G)]$, we have

$$\phi_{V,U}^G(v) = \sum_{(T,S) \in [\Pi(G)]} r_{T,S}\phi_{V,U}^G(\xi_{T,S}^G) = r_{V,U}.$$

It follows that for all $(V, U) \in [\Pi(G)]$

$$v.\xi_{V,U}^G = \phi_{V,U}^G(v)\xi_{V,U}^G.$$

Conversely, let $u$ be an element of $\mathbb{K}\Xi(G)$ satisfying $v.u = \lambda(v)u$ for any $v \in \mathbb{K}\Xi(G)$ where $\lambda(v) \in \mathbb{K}$. If $u = 0$, then there is nothing to prove. Otherwise, there exists a slice $(T, S)$ of $G$ such that $\xi_{T,S}^G.u \neq 0$. Then $\xi_{T,S}^G.u = \phi_{T,S}^G(u)\xi_{T,S}^G = \lambda(\xi_{T,S}^G)u \neq 0$, so $u$ is a scalar multiple of $\xi_{T,S}^G$. $\quad\square$

## 4. Effect of biset operations on idempotents

Recall that the algebra $\mathbb{K}\Xi(G)$ is commutative and split semi-simple, for any finite group $G$. Since $\mathbb{K}\Xi$ is a biset functor, by Remark 2.7, it becomes natural to look at the effect of elementary operations (Ind, Res, Inf, Def and Iso) on idempotents of the algebra $\mathbb{K}\Xi(G)$.

**Proposition 4.1.** *Let $(T, S)$ be a slice of $G$ and $H$ be subgroup of $G$. Then*

$$\mathrm{Res}_H^G \xi_{T,S}^G = \sum_{\substack{(T', S') \in [\Pi(H)] \\ (T', S') =_G (T, S)}} \xi_{T', S'}^H.$$

**Proof.** Let $(T, S)$ be slice of $G$.

By Proposition 3.15, the idempotent $\xi_{T,S}^G$ has the property that $\xi_{T,S}^G \pi(f) = \phi_{T,S}^G(\pi(f))\xi_{T,S}^G$ for any $\pi(f) \in \mathbb{K}\Xi(G)$. Observe that the restriction functor $\mathrm{Res}_H^G : G\text{-}\mathbf{set} \to H\text{-}\mathbf{set}$ induces a ring homomorphism $\mathrm{Res}_H^G : \mathbb{K}\Xi(G) \to \mathbb{K}\Xi(H)$. It follows that $\mathrm{Res}_H^G \xi_{T,S}^G$ is an idempotent of $\mathbb{K}\Xi(H)$, hence a sum of idempotents $\xi_{T', S'}^H$ for some slices $(T', S')$ of $H$. The idempotent $\xi_{V,U}^H$ appears in this sum if and only if

$$0 \neq \xi_{V,U}^H \mathrm{Res}_H^G \xi_{T,S}^G = \phi_{V,U}^H(\mathrm{Res}_H^G \xi_{T,S}^G)\xi_{V,U}^H = \phi_{V,U}^G(\xi_{T,S}^G)\xi_{V,U}^H$$

where the last equality follows from Proposition 3.6.

Thus $\xi_{V,U}^H$ appears in this sum if and only if $(V, U)$ is conjugate to $(T, S)$ in $G$. $\quad\square$

**Proposition 4.2.** *Let $H$ be a subgroup of $G$ and $(V, U)$ be a slice of $H$. Then*

$$\mathrm{Ind}_H^G \xi_{V,U}^H = \frac{|N_G(V, U)|}{|N_H(V, U)|} \xi_{V,U}^G.$$

**Proof.** It is straightforward from Notation 3.12 since for any element $\langle V, U \rangle_H$ of $\mathbb{K}\Xi(H)$, we have $\mathrm{Ind}_H^G(\langle V, U \rangle_H) = \langle V, U \rangle_G$. $\quad\square$

**Proposition 4.3.** *Let $N$ be a normal subgroup of $G$. Then for any slice $(T, S)$ of $G$ such that $N$ is contained in $S$, we have*

$$\mathrm{Inf}_{G/N}^G \xi_{T/N, S/N}^{G/N} = \sum_{\substack{(Y, X) \in [\Pi(G)] \\ (YN, XN) =_G (T, S)}} \xi_{Y,X}^G.$$

**Proof.** We observe that $\mathrm{Inf}_{G/N}^{G} : \mathbb{K}\Xi(G/N) \to \mathbb{K}\Xi(G)$ is a ring homomorphism.

It follows that $\mathrm{Inf}_{G/N}^{G}\xi_{T/N,S/N}^{G/N}$ is an idempotent of $\mathbb{K}\Xi(G)$, hence a sum of idempotents $\xi_{Y,X}^{G}$, for some slices $(Y,X) \in \Pi(G)$. The idempotent $\xi_{Y,X}^{G}$ appears in this sum if and only if

$$0 \neq \xi_{Y,X}^{G}\mathrm{Inf}_{G/N}^{G}\xi_{T/N,S/N}^{G/N} = \phi_{Y,X}^{G}(\mathrm{Inf}_{G/N}^{G}\xi_{T/N,S/N}^{G/N})\xi_{Y,X}^{G} = \phi_{YN/N,XN/N}^{G/N}(\xi_{T/N,S/N}^{G/N})\xi_{Y,X}^{G}$$

i.e. if and only if $(YN, XN)$ is conjugate to $(T,S)$ in $G$ i.e. if and only if $(Y,X) \twoheadrightarrow (T,S)$ modulo $G$. $\quad\square$

**Proposition 4.4.** *Let $N$ be a normal subgroup of $G$. Then*

$$\mathrm{Def}_{G/N}^{G}\xi_{G,S}^{G} = \mathrm{m}_{G,S,N}\xi_{G/N,SN/N}^{G/N}$$

*where*

$$\mathrm{m}_{G,S,N} = \frac{|N_G(SN):SN|}{|N_G(S)|} \sum_{\substack{U \leq S \leq V \leq G \\ VN=G \\ UN=SN}} |U|\mu(U,S)\mu(V,G).$$

*More generally, if $(T,S)$ is some slice of $G$, then*

$$\mathrm{Def}_{G/N}^{G}\xi_{T,S}^{G} = \frac{|N_T(S)||N_G(TN,SN)|}{|N_G(T,S)||N_T(SN)|}\mathrm{m}_{T,S,T\cap N}\xi_{TN/N,SN/N}^{G/N}.$$

**Proof.** We observe that if $X$ is a $G$-set, then $\mathrm{Def}_{G/N}^{G}X = (G/N) \times_G X \simeq N\backslash X$.

So the map $\mathrm{Def}_{G/N}^{G} : \mathbb{Q}\Xi(G) \to \mathbb{Q}\Xi(G/N)$ is such that $\mathrm{Def}_{G/N}^{G}\big(\pi(\ X \xrightarrow{f} Y\ )\big) = \pi(\ N\backslash X \xrightarrow{N\backslash f} N\backslash Y\ )$ for any morphism of $G$-sets $f : X \to Y$.

Now if $Z$ is a $G/N$-set, and $X$ is a $G$-set, then there is an isomorphism of $G/N$-sets

$$\alpha_{Z,X} : Z \times (N\backslash X) \to \big(N\backslash(\mathrm{Inf}_{G/N}^{G}Z) \times X)\big).$$

Let $Z \xrightarrow{g} T$ is a morphism of $G/N$-sets, and $X \xrightarrow{f} Y$ is a morphism of $G$-sets then the morphisms $g \times (N\backslash f)$ and $N\backslash\big((\mathrm{Inf}_{G/N}^{G}g) \times f\big)$ have the same image in $\mathbb{Q}\Xi(G/N)$. Then

$$v\mathrm{Def}_{G/N}^{G}u = \mathrm{Def}_{G/N}^{G}\big((\mathrm{Inf}_{G/N}^{G}v)u\big),$$

for any $v \in \mathbb{Q}\Xi(G/N)$ and $u \in \mathbb{Q}\Xi(G)$. For $u = \xi_{G,S}^{G}$, this gives

$$v\mathrm{Def}_{G/N}^G \xi_{G,S}^G = \mathrm{Def}_{G/N}^G \big( (\mathrm{Inf}_{G/N}^G v)\xi_{G,S}^G \big)$$

$$= \phi_{G,S}^G (\mathrm{Inf}_{G/N}^G v)\mathrm{Def}_{G/N}^G \xi_{G,S}^G$$

$$= \phi_{G/N,SN/N}^{G/N}(v)\mathrm{Def}_{G/N}^G \xi_{G,S}^G$$

and by Proposition 3.15 it follows that $\mathrm{Def}_{G/N}^G \xi_{G,S}^G$ is equal to a scalar multiple of $\xi_{G/N,SN/N}^{G/N}$. In particular, there is a scalar $m$ such that

$$\mathrm{Def}_{G/N}^G \xi_{G,S}^G = m\xi_{G/N,SN/N}^{G/N}$$

and the case $\mathbb{K} = \mathbb{Q}$ shows that $m \in \mathbb{Q}$. We set $\bar{G} = G/N$, $\bar{S} = SN/N$, $\bar{Y} = Y/N$ and $\bar{X} = X/N$.

Since $\mathrm{Def}_{G/N}^G \langle T, S \rangle_G = \langle TN, SN \rangle_G = \langle TN/N, SN/N \rangle_{G/N}$ for any slice $(T, S)$ of $G$, this yields

$$\text{(1)} \quad \mathrm{Def}_{G/N}^G \xi_{G,S}^G = \frac{1}{|N_G(S)|} \sum_{U \le S \le V \le G} |U|\mu(U,S)\mu(V,G)\langle VN/N, UN/N \rangle_{G/N}$$

$$= m\xi_{G/N,SN/N}^{G/N}$$

$$= \frac{m}{|N_{\bar{G}}(\bar{S})|} \sum_{\bar{X} \le \bar{S} \le \bar{Y} \le \bar{G}} |\bar{X}|\mu(\bar{X},\bar{S})\mu(\bar{Y},\bar{G})\langle \bar{Y}, \bar{X} \rangle_{\bar{G}}$$

$$\text{(2)} \qquad = \frac{m}{|N_{G/N}(SN/N)||N|} \sum_{\substack{X \le SN \le Y \le G \\ N \le X}} |X|\mu(X,SN)\mu(Y,G)\langle \bar{Y}, \bar{X} \rangle_{\bar{G}}.$$

The coefficient of $\langle G/N, SN/N \rangle_{G/N}$ is equal to $\frac{m|SN|}{|N_{G/N}(SN/N)||N|}$ in (2), and equal to

$$\frac{1}{|N_G(S)|} \sum_{\substack{U \le S \le V \le G \\ VN =_G \bar{G} \\ UN =_G SN}} |U|\mu(U,S)\mu(V,G)$$

in (1).

Since $U \le S$ and $UN =_G SN$ implies that $UN = SN$, and the coefficient of $\langle G/N, SN/N \rangle_{G/N}$ is equal to

$$\frac{1}{|N_G(S)|} \sum_{\substack{U \le S \le V \le G \\ VN = \bar{G} \\ UN = SN}} |U|\mu(U,S)\mu(V,G)$$

in (1).

It follows that

$$\mathrm{m}_{G,S,N} := m = \frac{|N_{G/N}(SN/N)||N|}{|N_G(S)||SN|} \sum_{\substack{U \le S \le V \le G \\ VN=G \\ UN=SN}} |U|\mu(U,S)\mu(V,G).$$

We observe also from (2) that,

$$\xi_{T,S}^G = \frac{|N_T(S)|}{|N_G(T,S)|}\mathrm{Ind}_T^G\xi_{T,S}^T,$$

for any slice $(T,S)$ of $G$. Thus, by [Lemma 2.5](#)

$$
\begin{aligned}
\mathrm{Def}_{G/N}^G\xi_{T,S}^G &= \frac{|N_T(S)|}{|N_G(T,S)|}\mathrm{Def}_{G/N}^G\mathrm{Ind}_T^G\xi_{T,S}^T \\
&= \frac{|N_T(S)|}{|N_G(T,S)|}\mathrm{Ind}_{TN/N}^{G/N}\mathrm{Iso}_{T/T\cap N}^{TN/N}\mathrm{Def}_{T/T\cap N}^T\xi_{T,S}^T \\
&= \frac{|N_T(S)|}{|N_G(T,S)|}\mathrm{m}_{T,S,T\cap N}\mathrm{Ind}_{TN/N}^{G/N}\mathrm{Iso}_{T/T\cap N}^{TN/N}\xi_{T/T\cap N,S.(T\cap N)/T\cap N}^{T/T\cap N} \\
&= \frac{|N_T(S)|}{|N_G(T,S)|}\mathrm{m}_{T,S,T\cap N}\mathrm{Ind}_{TN/N}^{G/N}\xi_{TN/N,SN/N}^{TN/N} \\
&= \frac{|N_T(S)||N_{G/N}(TN/N,SN/N)|}{|N_G(T,S)||N_{T/N}(SN/N)|}\mathrm{m}_{T,S,T\cap N}\xi_{TN/N,SN/N}^{G/N} \\
&= \frac{|N_T(S)||N_G(TN,SN)|}{|N_G(T,S)||N_T(SN)|}\mathrm{m}_{T,S,T\cap N}\xi_{TN/N,SN/N}^{G/N}.
\end{aligned}
$$

This completes the proof.    □

**Proposition 4.5.** *If $\varphi : G \to G'$ is a group isomorphism, and $(T,S)$ is a slice of $G$, then*

$$\mathrm{Iso}(\varphi)(\xi_{T,S}^G) = \xi_{\varphi(T),\varphi(S)}^{G'}.$$

**Proof.** This is straightforward.    □

## 5. Properties of $\mathrm{m}_{G,S,N}$'s

Let $(G,S)$ be a slice and $N$ be a normal subgroup of $G$.

**Proposition 5.1.** *Let $S \le G$ be finite groups. If $M$ and $N$ are normal subgroups of $G$ with $N \le M$, then*

$$\mathrm{m}_{G,S,M} = \mathrm{m}_{G,S,N}\mathrm{m}_{G/N,SN/N,M/N}.$$

**Proof.** This follows from the transitivity property of deflation maps:

$$\mathrm{Def}_{G/M}^{G}\xi_{G,S}^{G} = \mathrm{m}_{G,S,M}\xi_{G/M,SM/M}^{G/M}$$
$$= \mathrm{Def}_{G/M}^{G/N}\mathrm{Def}_{G/N}^{G}\xi_{G,S}^{G}$$
$$= \mathrm{m}_{G,S,N}\mathrm{Def}_{G/M}^{G/N}\xi_{G/N,SN/N}^{G/N}$$
$$= \mathrm{m}_{G,S,N}\mathrm{m}_{G/N,SN/N,M/N}\xi_{G/M,SM/M}^{G/M}. \qquad \square$$

**Remark 5.2.** If $G$ is a finite group, and $N$ is a normal subgroup of $G$. We observe that the constant $\mathrm{m}_{G,G,N}$ is equal to the constant $\mathrm{m}_{G,N}$ introduced in [1] and [4].

Recall some properties of this constant $\mathrm{m}_{G,N}$ (see [1]).

- For any finite group $G$ and any normal subgroup $N$ of $G$, we have $\mathrm{m}_{G,N} = \mathrm{m}_{G/\Phi(G),N\Phi(G)/\Phi(G)}$ where $\Phi(G)$ is the Frattini subgroup of $G$. So, if $G$ is a finite $p$-group (where $p$ is a prime), we can assume that $G$ is elementary abelian, in order to compute $\mathrm{m}_{G,N}$.
- If $G$ is an elementary abelian $p$-group of rank $n > 2$ and $N$ is a subgroup of $G$ of rank $k$ with $1 \leq k \leq n-2$, then

$$\mathrm{m}_{G,N} = \prod_{i=1}^{k}\left(1 - p^{n-1-i}\right).$$

In this case the constant $\mathrm{m}_{G,N}$ is equal to zero if and only if $G$ is noncyclic and the quotient $G/N$ is cyclic.

**Definition 5.3** *(see [2], Definition 5.4.6).* A finite group $G$ is called a **B**-group (over $\mathbb{K}$) if $|G| \neq 0$ in $\mathbb{K}$, and if for any non trivial normal subgroup $N$ of $G$, the constant $\mathrm{m}_{G,N}$ is equal to zero (in $\mathbb{K}$).

**Proposition 5.4** *(see [2], Section 5.6.9).* Let $G$ be a $p$-group. Then $G$ is a **B**-group if and only if $G$ is trivial or isomorphic to $C_p \times C_p$.

**Proposition 5.5.** Let $S \leq G$ be finite groups and $N$ be a normal subgroup of $G$. Then

$$\mathrm{m}_{G,S,N} = \frac{|N_G(SN) : SN|}{|N_G(S) : S|}\mathrm{m}_{S,S\cap N}\mathrm{m}_{G,S,N}^{\circ}$$

where $\mathrm{m}_{G,S,N}^{\circ} = \sum_{\substack{S \leq V \leq G \\ VN=G}} \mu(V,G).$

**Proof.** By Dedekind's identity we have,

$$\mathrm{m}_{S,S\cap N} = \frac{1}{|S|} \sum_{X(S\cap N)=S} |X|\mu(X,S)$$

$$= \frac{1}{|S|} \sum_{\substack{X\leq S \\ XN=SN}} |X|\mu(X,S).$$

We can now calculate the product

$$\mathrm{m}_{S,S\cap N}\mathrm{m}^{\circ}_{G,S,N} = \Big(\frac{1}{|S|} \sum_{X(S\cap N)=S} |X|\mu(X,S)\Big)\Big(\sum_{\substack{S\leq V\leq G \\ VN=G}} \mu(V,G)\Big)$$

$$= \Big(\frac{1}{|S|} \sum_{\substack{X\leq S \\ XN=SN}} |X|\mu(X,S)\Big)\Big(\sum_{\substack{S\leq V\leq G \\ VN=G}} \mu(V,G)\Big)$$

$$= \frac{1}{|S|} \sum_{\substack{X\leq S\leq V\leq G \\ VN=G \\ XN=SN}} |X|\mu(X,S)\mu(V,G)$$

$$= \frac{|N_G(S):S|}{|N_G(SN):SN|}\mathrm{m}_{G,S,N}. \qquad \square$$

**Proposition 5.6.** *Let $G$ be a finite group.*

*If $N$ is a minimal abelian normal subgroup of $G$, then*

$$\mathrm{m}_{G,1,N} = \frac{1}{|N|}\Big(1 - |K_G(N)|\Big) = \frac{1}{|N|}\mathrm{m}^{\circ}_{G,1,N},$$

*where $K_G(N)$ is the set of complements of $N$ in $G$.*

**Proof.** Let $X$ be a subgroup of $G$ such that $XN = G$. Then $X \cap N$ is normalized by $X$, and by $N$ since $N$ is abelian. It follows that $X \cap N \unlhd G$, thus $X \cap N = N$ or $X \cap N = 1$ by minimality of $N$. In the first case $X = G$, and in the second case $X \in K_G(N)$, so $|X| = |G:N|$, and moreover $X$ is a maximal subgroup of $G$, so $\mu_G(X,G) = -1$. This shows that

$$\mathrm{m}_{G,1,N} = \frac{1}{|N|}\Big(\sum_{XN=G} \mu_G(X,G)\Big) = \frac{1}{|N|}\Big(1 + \sum_{\substack{X,N=G \\ X\neq G}} \mu_G(X,G)\Big)$$

$$= \frac{1}{|N|}\Big(1 - |K_G(N)|\Big). \qquad \square$$

**Lemma 5.7.** *Let $S \leq G$ be finite groups and $N$ be a normal subgroup of $G$. Then*

$$\mathrm{m}^{\circ}_{G,S,N} = \mathrm{m}^{\circ}_{G/\Phi(G),S\Phi(G)/\Phi(G),N\Phi(G)/\Phi(G)}.$$

**Proof.** We set $\bar{G} = G/\Phi(G)$, $\bar{S} = S\Phi(G)/\Phi(G)$ and $\bar{N} = N\Phi(G)/\Phi(G)$. Since $\mu(V,G) = 0$ if $\Phi(G) \nleq V$, it follows that

$$
\begin{aligned}
\mathrm{m}_{G,S,N}^{\circ} &= \sum_{\substack{S \leq V \leq G \\ VN=G}} \mu(V,G) = \sum_{\substack{S\Phi(G) \leq V \leq G \\ VN=G}} \mu(V,G) \\
&= \sum_{\substack{S\Phi(G) \leq V \leq G \\ V.\Phi(G)N=G}} \mu(V,G) = \mathrm{m}_{G,S\Phi(G),N\Phi(G)}^{\circ} \\
&= \sum_{\substack{\bar{S} \leq \bar{V} \leq \bar{G} \\ \bar{V}.\bar{N}=\bar{G}}} \mu(\bar{V},\bar{G}) = \mathrm{m}_{\bar{G},\bar{S},\bar{N}}^{\circ}
\end{aligned}
$$

where $\bar{V}$ denotes $V/\Phi(G)$. $\quad\square$

## 6. The case of $p$-groups

Let $p$ be a prime number, and $G$ be a finite $p$-group. Let $N$ be a normal subgroup of $G$. Since,

$$
\mathrm{m}_{G,S,N}^{\circ} = \mathrm{m}_{G/\Phi(G),S\Phi(G)/\Phi(G),N\Phi(G)/\Phi(G)}^{\circ}
$$

for any subgroup $S$ of $G$, in order to compute $\mathrm{m}_{G,S,N}^{\circ}$ we may assume that $G$ is an elementary abelian $p$-group without lost of generality.

**Proposition 6.1.** *Let $G$ be an elementary abelian $p$-group of rank $n$ with $n \geq 1$, and let $N$ be subgroup of $G$ of rank $k$ with $k \geq 1$. Then*

$$
\mathrm{m}_{G,1,N}^{\circ} = \prod_{i=1}^{k} \left(1 - p^{n-i}\right).
$$

*More generally, if $S$ is a subgroup of $G$, then $\mathrm{m}_{G,S,N}^{\circ} = \mathrm{m}_{G/S,1,NS/S}^{\circ}$.*

**Proof.** We argue by induction on $|G|$. If $n = 1$ then by Proposition 5.6

$$
\mathrm{m}_{C_p,1,C_p}^{\circ} = 1 - 1 = 0.
$$

One can assume that the result holds for all finite groups $H$ with $|H| < |G|$.

We assume first that $k = 1$.

Among the $(p^n - 1)/(p-1)$ maximal subgroups of $G$, a total of $(p^{n-1}-1)/(p-1)$ contain $N$. Hence $N$ has $(p^n - p^{n-1})/(p-1) = p^{n-1}$ complements in $G$, and so Proposition 5.6 gives

$$
\mathrm{m}_{G,1,N}^{\circ} = 1 - p^{n-1} \quad \text{in this case.}
$$

We now assume $k > 1$.

Let $M \leq G$ with $|M| = p$. Then the induction hypothesis applied to the group $G/M$ implies that

$$\mathrm{m}^{\circ}_{G/M,1,N/M} = (1 - p^{n-2}) \ldots (1 - p^{n-k}).$$

And by Proposition 5.1

$$\mathrm{m}^{\circ}_{G,1,N} = \mathrm{m}^{\circ}_{G,1,M} \mathrm{m}^{\circ}_{G/M,1,N/M} = (1 - p^{n-1})(1 - p^{n-2}) \ldots (1 - p^{n-k}).$$

The second part is an immediate consequence of the definition of the Möbius function and the fact that $V \mapsto V/S$ is an isomorphism from the poset of subgroups $V$ of $G$ containing $S$ *such that* $VN = G$ to the poset of subgroups $V/S$ of $G/S$ *such that* $V/S.NS/S = G/S$. Thus

$$\mathrm{m}^{\circ}_{G,S,N} = \sum_{\substack{S \leq V \leq G \\ VN = G}} \mu(V, G) = \sum_{\substack{S \leq V \leq G \\ V\bar{S}N = G}} \mu(V, G)$$

$$= \sum_{\substack{1 \leq V/S \leq G/S \\ V/\bar{S}.SN/\bar{S} = G/S}} \mu(V/S, G/S) = \mathrm{m}^{\circ}_{G/S,1,NS/S}. \qquad \square$$

**Remark 6.2.** The previous results show that if $G$ is a $p$-group and $S$, $N$ are subgroups of $G$, then: the constant $\mathrm{m}^{\circ}_{G,S,N}$ is equal to zero if and only if $G \neq S$ and $G = SN$.

**Remark 6.3.** Let $T$ be a finite $p$-group and let $S$ be a subgroup of $T$. Then for any normal subgroup $N$ of $T$,

$$\mathrm{m}_{T,S,N} = 0 \Leftrightarrow \begin{cases} S \text{ noncyclic} \\ \\ S/S \cap N \text{ cyclic} \end{cases} \text{ or } \begin{cases} T \neq S \\ \\ T = SN. \end{cases}$$

## 7. Characterization of the ideals of $\mathbb{K}\Xi$

In this section, we aim to give a combinatorial description of the ideals of the Green slice Burnside functor $\mathbb{K}\Xi$.

First we recall from Section 3 that for any finite group $G$, the algebra $\mathbb{K}\Xi(G)$ is a split semisimple commutative $\mathbb{K}$-algebra, whose primitive idempotents $\xi^G_{T,S}$ are indexed by a set of representatives of conjugacy classes of slices of $G$. In particular, if $F$ is an ideal of $\mathbb{K}\Xi$, we denote by $\mathcal{A}_{F,G}$ the set of slices $(T, S)$ of $G$ such that $\xi^G_{T,S} \in F(G)$. Then $F(G)$ is the sum of one-dimensional $\mathbb{K}$-vector spaces generated by $\xi^G_{T,S}$, for $(T, S)$ in the family $\mathcal{A}_{F,G}$. This sum is direct if the slices $(T, S)$ are chosen up to conjugation in $G$.

Our next result will be essential in our analysis of the ideals of $\mathbb{K}\Xi$.

**Theorem 7.1.** *If $F$ is an ideal of the Green functor $\mathbb{K}\Xi$, then $\mathcal{A}_{F,G}$ satisfies the four following conditions:*

**A** Let $(T, S)$ be an element of $\mathcal{A}_{F,G}$ and $(T', S')$ be a slice of $G$. If $(T, S)$ and $(T', S')$ are isomorphic, then $(T', S')$ is also an element of $\mathcal{A}_{F,G}$.

**B** The family $\mathcal{A}_{F,G}$ depends only on $F$ and not on the group $G$. So it will be denoted by $\mathcal{A}_F$ and $\mathcal{A}_F = \{(T, S) \mid \xi_{T,S}^T \in F(T)\}$.

**C** Let $N$ be a normal subgroup of $T$ and let $(T/N, S/N)$ be an element of $\mathcal{A}_F$. If there exist a slice $(Y, X)$ and a surjective group homomorphism $(Y, X) \twoheadrightarrow (T, S)$, then $(Y, X)$ lies in $\mathcal{A}_F$.

**D** Let $N$ be a normal subgroup of $T$ and let $(T, S)$ be an element of $\mathcal{A}_F$ such that the constant $\mathrm{m}_{T,S,N}$ is different from zero. Then the slice $(T/N, SN/N)$ lies in $\mathcal{A}_F$.

*Conversely, if a family $\mathcal{A}$ of slices is given, we set*

$$F_\mathcal{A}(G) = \sum_{(T,S) \in \mathcal{A}} \mathbb{K} \xi_{T,S}^G,$$

*for any finite group $G$. If $\mathcal{A}$ fulfills Conditions **A** to **D**, then $F_\mathcal{A}$ is an ideal of the Green functor $\mathbb{K}\Xi$.*

**Proof.** Let $F$ be an ideal of $\mathbb{K}\Xi$. We claim that Conditions **A** to **D** are fulfilled.

Condition **A** is obvious by Proposition 4.5.

For each $T \le G$, we must have $\mathrm{Ind}_T^G F(T) \subseteq F(G)$. Let $(T, S)$ be an element of $\mathcal{A}_{F,T}$. By Proposition 4.2, the idempotent $\xi_{T,S}^G$ is a non zero scalar multiple of $\mathrm{Ind}_T^G \xi_{T,S}^T$. Hence $(T, S) \in \mathcal{A}_{F,G}$.

We must also have $\mathrm{Res}_T^G F(G) \subseteq F(T)$. Let $(T, S)$ be an element of $\mathcal{A}_{F,G}$.

Then by Proposition 4.1, the idempotent $\mathrm{Res}_T^G \xi_{T,S}^G$ is a sum of idempotents $\xi_{T',S'}^T$ (with non zero coefficient), for slices (T',S') of $T$ such that $(T', S') =_G (T, S)$.

Since $\xi_{T,S}^T$ appears in this sum with non zero coefficient, we have that $(T, S) \in \mathcal{A}_{F,T}$.

Thus, the slice $(T, S) \in \mathcal{A}_{F,G}$ if and only if it lies in $\mathcal{A}_{F,T}$. In other words the family $\mathcal{A}_{F,G}$ depends only on $F$ and $\mathcal{A}_{F,G} = \mathcal{A}_F = \{(T, S) \mid \xi_{T,S}^T \in F(T)\}$.

Condition **C** and Condition **D** follow obviously from Proposition 4.3 and Proposition 4.4.

Conversely, let $\mathcal{A}$ be a set of slices which satisfies Conditions **A**, **B**, **C** and **D**. In order to prove that $F_\mathcal{A}$ is an ideal of $\mathbb{K}\Xi$, all we have to show is that it is a biset subfunctor of $\mathbb{K}\Xi$, since for any $G$, the vector space $F_\mathcal{A}(G)$ is an ideal of $\mathbb{K}\Xi(G)$.

As the elementary bisets generate the biset category, it is enough to check that $F_\mathcal{A}$ is stable under transport by isomorphism, induction, restriction, inflation and deflation.

The first case is clear.

Let $H$ be a subgroup of $G$ and $(V, U)$ be a slice of $H$ such that $\xi_{V,U}^H \in F_\mathcal{A}(H)$. Then $\xi_{V,U}^V \in F_\mathcal{A}(V)$ and $\xi_{V,U}^G \in F_\mathcal{A}(G)$, by Condition **B**.

Since, the idempotent $\xi_{V,U}^G$ is a non zero scalar multiple of $\mathrm{Ind}_H^G \xi_{V,U}^V$, we observe that $\mathrm{Ind}_H^G \xi_{V,U}^H \in F_\mathcal{A}(G)$ and so $\mathrm{Ind}_H^G F_\mathcal{A}(H) \subseteq F_\mathcal{A}(G)$.

Now, let $(T, S)$ be a slice of $G$ such that $\xi_{T,S}^G \in F_\mathcal{A}(G)$ and $H$ be a subgroup of $G$. By the Mackey formula the restriction of $\xi_{T,S}^G$ to $H$ is equal to

$$\operatorname{Res}_H^G \xi_{T,S}^G = \frac{|N_H(T,S)|}{|N_G(T,S)|} \operatorname{Res}_H^G \operatorname{Ind}_T^G \xi_{T,S}^T$$

$$= \frac{|N_H(T,S)|}{|N_G(T,S)|} \sum_{x \in [H \backslash G/T]} \operatorname{Ind}_{H \cap {}^x T}^H \circ \operatorname{Iso}(\gamma_x) \circ \operatorname{Res}_{H^x \cap T}^T \xi_{T,S}^T$$

where $[H \backslash G/T]$ is a set of representatives of $(H,T)$-double cosets in $G$, and $\gamma_x$ is the group isomorphism induced by conjugation by $x$.

We observe by Condition **B**, that $\xi_{T,S}^T \in F_{\mathcal{A}}(T)$. Since the restriction $\operatorname{Res}_{H^x \cap T}^T \xi_{T,S}^T$ is either zero or $\xi_{T,S}^T$, then $\operatorname{Res}_{H^x \cap T}^T \xi_{T,S}^T$ lies in $F_{\mathcal{A}}(T)$.

By Condition **A**, the element $\operatorname{Iso}(\gamma_x) \circ \operatorname{Res}_{H^x \cap T}^T \xi_{T,S}^T$ lies in $F_{\mathcal{A}}(H \cap {}^x T)$ and from the above proof, the element

$$\operatorname{Ind}_{H \cap {}^x T}^H \circ \operatorname{Iso}(\gamma_x) \circ \operatorname{Res}_{H^x \cap T}^T \xi_{T,S}^T$$

is an element of $F_{\mathcal{A}}(H)$.

Thus $\operatorname{Res}_H^G F_{\mathcal{A}}(G) \subseteq F_{\mathcal{A}}(H)$, for any subgroup $H$ of $G$.

Next we show that $F_{\mathcal{A}}$ is stable under inflation. Let $N$ be a normal subgroup of $G$ and let $(T,S)$ be a slice of $G$ such that $N$ is contained in $S$. By Proposition 4.3, we have

$$\operatorname{Inf}_{G/N}^G \xi_{T/N,S/N}^{G/N} = \sum_{\substack{(Y,X) \in [\Pi(G)] \\ (Y,X) \twoheadrightarrow (T,S)}} \xi_{Y,X}^G.$$

Assume that $\xi_{T/N,S/N}^{G/N} \in F_{\mathcal{A}}(G/N)$. Then by Condition **C**, for any slice $(Y,X)$ of $G$ such that $(Y,X) \twoheadrightarrow (T,S)$, we obtain $\xi_{Y,X}^Y \in F_{\mathcal{A}}(Y)$. And by Condition **B**, we have $\xi_{Y,X}^G \in F_{\mathcal{A}}(G)$. Thus, all these $\xi_{Y,X}^G$ lie in $F_{\mathcal{A}}(G)$. Therefore $\operatorname{Inf}_{G/N}^G F_{\mathcal{A}}(G/N) \subseteq F_{\mathcal{A}}(G)$.

Now we prove that $F_{\mathcal{A}}$ is stable under deflation. Let $N$ be a normal subgroup of $G$ and let $(T,S)$ be a slice of $G$. Since

$$\operatorname{Def}_{G/N}^G \xi_{T,S}^G = \frac{|N_T(S)|}{|N_G(T,S)|} \operatorname{Ind}_{TN/N}^{G/N} \operatorname{Iso}_{T/T \cap N}^{TN/N} \operatorname{Def}_{T/T \cap N}^T \xi_{T,S}^T,$$

it suffices to show that $\operatorname{Def}_{T/T \cap N}^T F_{\mathcal{A}}(T) \subseteq F_{\mathcal{A}}(T/T \cap N)$. Recall from Proposition 4.4 that

$$\operatorname{Def}_{T/T \cap N}^T \xi_{T,S}^T = \mathrm{m}_{T,S,T \cap N} \xi_{T/T \cap N, S(T \cap N)/T \cap N}^{T/T \cap N}.$$

Thus, if $\xi_{T,S}^T \in F_{\mathcal{A}}(T)$, then by Condition **D** the idempotent $\xi_{T/T \cap N, S(T \cap N)/T \cap N}^{T/T \cap N} \in F_{\mathcal{A}}(T/T \cap N)$ whenever $\mathrm{m}_{T,S,T \cap N} \neq 0$.

Thus $\operatorname{Def}_{T/T \cap N}^T F_{\mathcal{A}}(T) \subseteq F_{\mathcal{A}}(T/T \cap N)$. This completes the proof. $\quad \square$

**Notation 7.2.** For any ideal $F$ of $\mathbb{K}\Xi$, we denote by $\mathcal{A}_F$ the family of slices corresponding to $F$ i.e. such that

$$F(G) = \sum_{(T,S) \in \mathcal{A}_F} \mathbb{K} \xi_{T,S}^G$$

for any finite group $G$.

**Corollary 7.3.** *Let $F$ be an ideal of $\mathbb{K}\Xi$.*

*If $(T, S)$ is an element of $\mathcal{A}_F$ and $(B, A)$ is any slice, then the slice $(B \times T, A \times S)$ lies in $\mathcal{A}_F$.*

**Proof.** The subgroup $N = B \times 1$ of $B \times T$ is normal and,

$$(B \times T)/(B \times 1) \simeq T \quad \text{and} \quad \big((A \times S)(B \times 1)\big)/(B \times 1) \simeq S.$$

Thus $(B \times T, A \times S) \twoheadrightarrow (T, S)$, and by Condition **C** we have $(B \times T, A \times S) \in \mathcal{A}_F$. $\quad\square$

**Corollary 7.4.** *Let $F$ be an ideal of $\mathbb{K}\Xi$.*

*If $T$ is minimal such that there exists $S$ with $(T, S) \in \mathcal{A}_F$, then for any non trivial normal subgroup $N$ of $T$, the constant $\mathrm{m}_{T,S,N}$ is equal to zero.*

**Proof.** This follows from Condition **D** because, if $1 \neq N \trianglelefteq T$ and $\mathrm{m}_{T,S,N} \neq 0$ then $(T/N, SN/N) \in \mathcal{A}_F$. $\quad\square$

The above corollary motivates the following definition:

**Definition 7.5.** A slice $(T, S)$ of $T$ is called a **T**-slice (over $\mathbb{K}$) if for any non trivial normal subgroup $N$ of $T$, the constant $\mathrm{m}_{T,S,N}$ is equal to zero.

## 8. Ideals of the slice Burnside $p$-biset functor $\mathbb{K}\Xi_p$

Throughout this section $p$ denotes a prime number.

We apply the results of the previous sections to the determination of the full lattice of ideals of the slice Burnside $p$-biset functor. Unless otherwise specified, the groups considered in this section are $p$-groups.

**Definition 8.1.**

- The $p$-biset category $\mathbb{K}\mathcal{C}_p$ is the full subcategory of the biset category $\mathbb{K}\mathcal{C}$ whose objects are finite $p$-groups.
- A $p$-biset functor on $\mathcal{C}_p$ with values in $\mathbb{K}$-Vect is a $\mathbb{K}$-linear functor from $\mathbb{K}\mathcal{C}_p$ to $\mathbb{K}$-Vect.

**Remark 8.2.** By restricting to $p$-groups only, the category of bisets functors becomes the category of $p$-biset functors. All the methods of biset functors can be easily adapted and therefore our results hold in this context. In particular:

- The $p$-biset functors form an abelian category.
- The correspondence $G \mapsto \mathbb{K}\Xi(G)$ is a biset functor on the full subcategory of the biset category consisting of finite $p$-groups. It is denoted $\mathbb{K}\Xi_p$ and is called slice Burnside $p$-biset functor.
- As before, for any ideal $F$ of $\mathbb{K}\Xi_p$, we denote by $\mathcal{A}_F$ the family of slices corresponding to $F$.

**Proposition 8.3.** *The following families of slices*

1. $\mathcal{A}_{\mathfrak{I}_1} = \{(T, S) \mid T \neq S\}$,
2. $\mathcal{A}_{\mathfrak{I}_2} = \{(T, S) \mid S \ noncyclic\}$,
3. $\mathcal{A}_{\mathfrak{I}_3} = \{(T, S) \mid T \neq S, \ S \ noncyclic\}$

*satisfy conditions* **A**, **B**, **C** *and* **D**.

**Proof.** It is clear for **A**, **B** and **C**.

In order to prove **D** for the families $\mathcal{A}_{\mathfrak{I}_1}$, $\mathcal{A}_{\mathfrak{I}_2}$ and $\mathcal{A}_{\mathfrak{I}_3}$, let $(T, S)$, $(V, U)$ and $(B, A)$ be slices such that $(T, S) \in \mathcal{A}_{\mathfrak{I}_1}$, $(V, U) \in \mathcal{A}_{\mathfrak{I}_2}$ and $(B, A) \in \mathcal{A}_{\mathfrak{I}_3}$.

1. Then, if $N$ is normal subgroup of $T$ such that the constant $\mathrm{m}_{T,S,N}$ is nonzero then

$$\mathrm{m}_{T,S,N}^{\circ} = \mathrm{m}_{T/\Phi(T),S\Phi(T)/\Phi(T),N\Phi(T)/\Phi(T)}^{\circ} \neq 0.$$

By Proposition 6.1, this is equivalent to saying that either

$$T/\Phi(T) = \big(S\Phi(T)\big)/\Phi(T) \ \text{ or } \ \Big(\big(S\Phi(T)\big)/\Phi(T)\Big).\Big(\big(N\Phi(T)\big)/\Phi(T)\Big) \neq T/\Phi(T).$$

Since $T \neq S$, we have $S\Phi(T) \neq T$. It follows from the previous argument that $SN\Phi(T) \neq T$, and so $SN/N \neq T/N$. Thus $(T/N, SN/N) \in \mathcal{A}_{\mathfrak{I}_1}$.

2. If $N$ is normal subgroup of $V$ such that $\mathrm{m}_{V,U,N} \neq 0$ then $\mathrm{m}_{U,UN/N} = \mathrm{m}_{\bar{U},\bar{U}\bar{N}/\bar{N}} \neq 0$ where we denote $\bar{U} = U\Phi(T)/\Phi(T)$ and $\bar{N} = N\Phi(T)/\Phi(T)$. Since $U$ is noncyclic, the quotient $\bar{U}$ is noncyclic and by Proposition 6.1, we have $\bar{U}\bar{N}/\bar{N}$ noncyclic. Thus $UN/N$ is noncyclic and $(V/N, UN/N) \in \mathcal{A}_{\mathfrak{I}_2}$.

3. To prove **D** for $\mathfrak{I}_3$, it suffices to observe that $\mathcal{A}_{\mathfrak{I}_3} = \mathcal{A}_{\mathfrak{I}_1} \cap \mathcal{A}_{\mathfrak{I}_2}$. $\quad\square$

**Remark 8.4.** Note that by Theorem 7.1, the correspondence sending a finite $p$-group $G$ to the $\mathbb{K}$-vector space $\mathfrak{I}_1(G) = \sum_{(T,S)\in\mathcal{A}_{\mathfrak{I}_1}} \mathbb{K}\xi_{T,S}^G$ (resp. $\mathfrak{I}_2(G) = \sum_{(T,S)\in\mathcal{A}_{\mathfrak{I}_2}} \mathbb{K}\xi_{T,S}^G$ or $\mathfrak{I}_3(G) = \sum_{(T,S)\in\mathcal{A}_{\mathfrak{I}_3}} \mathbb{K}\xi_{T,S}^G$) is an ideal of $\mathbb{K}\Xi_p$.

**Proposition 8.5.** *Let $F$ be an ideal of the Green $p$-biset functor $\mathbb{K}\Xi_p$ and $(T, S)$ be a slice. The family $\mathcal{A}_F$ has the following property:*

1. *If $S$ is cyclic then*

$$(T, S) \in \mathcal{A}_F \Leftrightarrow \big(T/\Phi(T), S\Phi(T)/\Phi(T)\big) \in \mathcal{A}_F.$$

2. *If $E$ is an elementary abelian p-group, and if $S$ and $E \times S\Phi(T)/\Phi(T)$ are noncyclic, then*

$$(T, S) \in \mathcal{A}_F \Leftrightarrow \Big( E \times \big( T/\Phi(T) \big), E \times \big( S\Phi(T)/\Phi(T) \big) \Big) \in \mathcal{A}_F.$$

**Proof.** Let us point out that the first assertion holds for arbitrary finite groups.

1. Assume that $S$ is cyclic and the slice $(T, S) \in \mathcal{A}_F$.

   Since $S$ is cyclic, the constant $\mathrm{m}_{S, S \cap N}$ is different from zero, for any normal subgroup $N$ of $T$. Now by Proposition 5.5, for any normal subgroup $N$ of $T$, we have

   $$\mathrm{m}_{T, S, N} = \frac{|N_T(SN) : SN|}{|N_T(S) : S|} \mathrm{m}_{S, S \cap N} \mathrm{m}_{T, S, N}^{\circ}.$$

   In particular, if $N$ is equal to the Frattini subgroup $\Phi(T)$ of $T$, then $N$ is normal and the constant

   $$\mathrm{m}_{T, S, N}^{\circ} = \sum_{\substack{S \leq V \leq T \\ VN = T}} \mu(V, T) = \mu(T, T) = 1.$$

   Hence $\mathrm{m}_{T, S, N}$ is non-zero and by Condition **D**, the slice $\Big( T/\Phi(T), S\Phi(T)/\Phi(T) \Big)$ lies in $\mathcal{A}_F$.

   Conversely, assume that the slice $\big( T/\Phi(T), S\Phi(T)/\Phi(T) \big)$ lies in $\mathcal{A}_F$. Then, since there is a surjective group homomorphism $(T, S) \twoheadrightarrow \big( T/\Phi(T), S\Phi(T)/\Phi(T) \big)$, Condition **C** implies that the slice $(T, S)$ lies in $\mathcal{A}_F$.

2. Let $E$ be an elementary abelian group and $(T, S)$ be a slice such that the groups $S$ and $E \times \big( S\Phi(T)/\Phi(T) \big)$ are noncyclic. Assume that the slice $(T, S)$ lies in $\mathcal{A}_F$. Then by Corollary 7.3 the slice $(E \times T, E \times S)$ is an element of $\mathcal{A}_F$. Let $N = \Phi(E \times T) = 1 \times \Phi(T)$. Then by definition

   $$\mathrm{m}_{E \times T, E \times S, N}^{\circ} = \sum_{\substack{E \times S \leq L \leq E \times T \\ L(1 \times \Phi(T)) = E \times T}} \mu(L, E \times T) = \mu(E \times T, E \times T) = 1.$$

   Since the groups

   $$E \times S \qquad \text{and} \qquad (E \times S)/(1 \times S \cap \Phi(T)) \simeq E \times \big( S\Phi(T)/\Phi(T) \big)$$

   are not cyclic, the constant $\mathrm{m}_{E \times S, 1 \times S \cap \Phi(T)}$ is nonzero. Hence by Proposition 5.5, the constant $\mathrm{m}_{E \times T, E \times S, N}$ is nonzero. Thus by Condition **D** the slice

   $$\Big( E \times \big( T/\Phi(T) \big), E \times \big( S\Phi(T)/\Phi(T) \big) \Big)$$

   lies in $\mathcal{A}_F$.

Conversely, assume that $\left(E \times \left(T/\Phi(T)\right), E \times \left(S\Phi(T)/\Phi(T)\right)\right) \in \mathcal{A}_F$. Then by Condition **C**, we have $(E \times T, E \times S) \in \mathcal{A}_F$. Since $(E \times S)/(E \times 1) \simeq S$ is noncyclic and

$$\mathrm{m}^{\circ}_{E \times T, E \times S, E \times 1} = \sum_{\substack{E \times S \leq L \leq E \times T \\ L(E \times 1) = E \times T}} \mu(L, E \times T) = \mu(E \times T, E \times T) = 1,$$

the constant $\mathrm{m}_{E \times T, E \times S, E \times 1}$ is nonzero. Thus, Condition **D** implies that the slice $(T, S)$ lies in $\mathcal{A}_F$, and this completes the proof. $\quad\square$

This Proposition shows that we can work with elementary abelian $p$-groups in order to describe all the ideals of $\mathbb{K}\Xi_p$.

**Theorem 8.6.** *Let $E$ be an elementary abelian $p$-group of rank $n$, and let $F$ be a subgroup of $E$. Then the slice $(E, F)$ is a **T**-slice if and only if up to isomorphism, it belongs to the set*

$$\{(1,1); (C_p, 1); (C_p^2, C_p^2); (C_p^3, C_p^2)\}.$$

**Proof.** Let $(E, F)$ be a **T**-slice.

- If $F$ is cyclic and $F \neq E$.
  Then, for any non-trivial subgroup $N$ of $E$, we have $\mathrm{m}_{E,F,N} = 0$. Since $F$ is cyclic by Proposition 5.5, this is equivalent to saying that $FN = E$, for any subgroup $N$ such that $|N| = p$.
  If the order of $F$ were $p^k$ with $k \geq 1$, we would have a subgroup $M$ of $F$ of order $p$, so that $M$ would verify $FM = F = E$. This cannot occur, since by hypothesis $F \neq E$. Thus $F$ is the trivial group. Since $FN = E$, for any subgroup $N$ of order $p$, so the only possibility for $E$ is $E \simeq C_p$ i.e.

$$(E, F) \simeq (C_p, 1).$$

- If $F$ is noncyclic and $F \neq E$.
  Since $(E, F)$ is a **T**-slice,
  (∗) for any subgroup N such that $|N| = p$, either $F/(F \cap N)$ is cyclic or $FN = E$.
  – The first step is to consider the possibility that the quotient $F/(F \cap N)$ is noncyclic. Then $FN = F \times N = E$. Since $F$ is not cyclic, there exists an integer $k > 1$ such that $F \simeq C_p^k$.
  Thus, we may find a subgroup $M$ of $F$ of order $p$, and by (∗)

$$FM = E \text{ or } F/M \text{ is cyclic.}$$

  But, since $FM = F \neq E$, we have that $F/M$ is cyclic. Thus, there exists $n \geq 1$ such that

$$F \simeq (F/M) \times M \simeq C_{p^n} \times C_p.$$

Since E is elementary abelian, we deduce that

$$F \simeq C_p \times C_p \quad \text{and} \quad E \simeq C_p \times C_p \times C_p.$$

Hence $(E, F) \simeq (C_p^3, C_p^2)$.

– Now consider the possibility that the quotient $F/(F \cap N) = F/N$ is cyclic. Then $F \cap N$ is non-trivial and $N \subseteq F$. Therefore, the quotient $F/(F \cap N)$ is cyclic. Then $F$ is isomorphic to $C_{p^n} \times C_p$, for some integer $n \geq 1$. Since E is elementary abelian, we have

$$F \simeq C_p \times C_p.$$

On the other hand, since $F \neq E$ and $E$ is elementary abelian, there exists a subgroup $M$ of $E$ of order $p$ such that $M \nsubseteq F$. By the hypothesis $(*)$, the quotient $F/(F \cap M)$ is cyclic or $FM = E$. But, the quotient $F/F \cap M$ is noncyclic because $F \cap M$ is trivial and $F$ is not cyclic. Thus

$$E \simeq C_p \times C_p \times C_p.$$

Hence $(E, F) \simeq (C_p^3, C_p^2)$.

• If $E = F \neq 1$ then, for any subgroup $N$ of $E$, we have $m_{E,E,N}^\circ = 1$. It follows that the slice $(E, E)$ is a **T**-slice if and only if $E$ is a **B**-group. It other words $(E, E)$ is a **T**-slice if and only if $E \simeq C_p \times C_p$ or $E = 1$. Thus

$$(E, E) \simeq (C_p \times C_p, C_p \times C_p).$$

Conversely, we shall show that the slices $(1, 1)$, $(C_p, 1)$, $(C_p^2, C_p^2)$ and $(C_p^3, C_p^2)$ are **T**-slices.

This follows directly from the expression of $m_{G,S,N}$ as a product of a nonzero constant by $m_{S,S \cap N} m_{G,S,N}^\circ$, for any slice $(G, S)$ and any normal subgroup $N$ of $G$.

By Proposition 6.1, we have

$$m_{C_p,1,C_p} = m_{C_p,1,C_p}^\circ = 0.$$

Since $C_p$ is the only non trivial normal subgroup of $C_p$, we claim that the slice $(C_p, 1)$ is a **T**-slice.

Let $N$ be a non trivial normal subgroup of $C_p^2$, then

$$m_{C_p^2,N} = m_{C_p^2,N \cap C_p^2} = 0 \quad \text{and} \quad m_{C_p^2,C_p^2,N} = 0.$$

Let $M$ be a non trivial normal subgroup of $C_p^3$, then $m_{C_p^3,C_p^2,M}$ has the following decomposition:

$$\mathrm{m}_{C_p^3, C_p^2, M} = (*)\mathrm{m}_{C_p^2, M \cap C_p^2} \mathrm{m}_{C_p^3, C_p^2, M}^\circ,$$

where $(*)$ is a non zero constant.

If $M \cap C_p^2 \neq 1$, then $\mathrm{m}_{C_p^2, M \cap C_p^2}$ is equal to zero by Proposition 5.4. It follows that the constant $\mathrm{m}_{C_p^3, C_p^2, M}$ is equal to zero.

If $M \cap C_p^2 = 1$, then $C_p^3 = C_p^2 \times M$ and by Proposition 6.1, we deduce that

$$\mathrm{m}_{C_p^3, C_p^2, M}^\circ = 0 \quad \text{and} \quad \mathrm{m}_{C_p^3, C_p^2, M} = 0.$$

This shows that the slice $(C_p^3, C_p^2)$ is a **T**-slice. $\quad\square$

**Notation 8.7.** Let $G$ be a $p$-group and let $(T, S)$ be a slice of $G$.

For any set $I$, let $(F_\alpha)_{\alpha \in I}$ be a family of ideal of $\mathbb{K}\Xi_p$ such that $F_\alpha(T) \ni \xi_{T,S}^T$. Then $\bigcap_{\alpha \in I} F_\alpha$ is an ideal of $\mathbb{K}\Xi_p$ and

$$\Big( \bigcap_{\alpha \in I} F_\alpha \Big)(T) = \bigcap_{\alpha \in I} \big( F_\alpha(T) \big) \ni \xi_{T,S}^T.$$

The ideal $\mathbb{E}_{T,S}$ of the Green functor $\mathbb{K}\Xi_p$ generated by $\xi_{T,S}^T$ is by definition the intersection of all ideals $F_\alpha$ of $F$ such that $\xi_{T,S}^T \in F(T)$.

**Proposition 8.8.** *Let $F$ be an ideal of $\mathbb{K}\Xi$. Then*

$$F = \sum_{\substack{(T,S) \ slices \\ \xi_{T,S}^T \in F(T)}} \mathbb{E}_{T,S}.$$

**Proof.** For any finite group $G$, the ideal $F(G)$ is the sum in $\mathbb{K}\Xi(G)$ of those one dimensional subspaces $\mathbb{K}\xi_{T,S}^G$ for which $\xi_{T,S}^G \in F(G)$. We observe by Theorem 7.1 that

$$\xi_{T,S}^G \in F(G) \Leftrightarrow \xi_{T,S}^T \in F(T).$$

If $\xi_{T,S}^G \in F(G)$, then $\mathrm{Res}_T^G \xi_{T,S}^G = \sum_{\substack{S' \bmod T \\ S' =_G S}} \xi_{T,S'}^T \in F(T)$ i.e. $\mathbb{E}_{T,S} \subseteq F$. Conversely, if $\mathbb{E}_{T,S} \subseteq F$, i.e. if $\xi_{T,S}^T \in F(T)$, then $\xi_{T,S}^G \in F(G)$. It follows that

$$F(G) = \sum_{\substack{(T,S) \in \Pi(G) \\ \xi_{T,S}^T \in F(T)}} \mathbb{E}_{T,S}(G),$$

for any finite group $G$. Thus

$$F = \sum_{\substack{(T,S) \ slices \\ \xi_{T,S}^T \in F(T)}} \mathbb{E}_{T,S}. \quad\square$$

**Theorem 8.9.** *The full lattice of ideals of $\mathbb{K}\Xi_p$ has the following structure:*

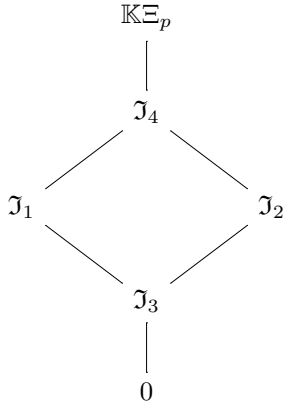

*where*

$$\mathbb{K}\Xi_p = \mathbb{E}_{1,1},$$

$$\mathfrak{I}_1 = \mathbb{E}_{C_p,1},$$

$$\mathfrak{I}_2 = \mathbb{E}_{C_p^2,C_p^2},$$

$$\mathfrak{I}_3 = \mathbb{E}_{C_p^3,C_p^2} = \mathfrak{I}_1 \cap \mathfrak{I}_2 \ and \ \mathfrak{I}_4 = \mathfrak{I}_1 + \mathfrak{I}_2.$$

*Moreover, the following equalities hold:*

1. $\mathcal{A}_{\mathbb{E}_{C_p,1}} = \{(T,S) \mid T \neq S\}$.
2. $\mathcal{A}_{\mathbb{E}_{C_p^2,C_p^2}} = \{(T,S) \mid S \ noncyclic\}$.
3. $\mathcal{A}_{\mathbb{E}_{C_p^3,C_p^2}} = \{(T,S) \mid T \neq S, \ S \ noncyclic\}$.

**Proof.** By Proposition 8.8, any ideal $F$ of $\mathbb{K}\Xi_p$ is equal to the sum of the ideals $\mathbb{E}_{T,S}$ it contains. Now by Proposition 8.5, for any slice $(T,S)$, there is a slice $(E,F)$, where $E$ is elementary abelian, such that $\mathbb{E}_{T,S} = \mathbb{E}_{E,F}$. Moreover, by Conditions **C** and **D**, we can assume that $(E,F)$ is a **T**-slice, i.e. by Theorem 8.6, one of the slices $(1,1)$, $(C_p,1)$, $(C_p^2,C_p^2)$, $(C_p^3,C_p^2)$.

- The first part of this Theorem follows easily from Conditions **A**, **B**, **C**, **D** and Corollary 7.3.

  Indeed, for any slice $(V,U)$, we have the projection $(V,U) \twoheadrightarrow (V/V, UV/V) = (1,1)$. Since $(1,1) \in \mathcal{A}_{\mathbb{E}_{1,1}}$ Condition **C** implies that $(V,U)$ lies in $\mathcal{A}_{\mathbb{E}_{1,1}}$. Thus $\mathbb{E}_{1,1} = \mathbb{K}\Xi_p$. Using Corollary 7.3 and the decomposition $(C_p^3,C_p^2) = (C_p^2,C_p^2) \times (C_p,1)$ it is easy to see that

$$\mathfrak{I}_3 \subseteq \mathfrak{I}_2 \ \text{and} \ \mathfrak{I}_3 \subseteq \mathfrak{I}_1$$

  as ideals.

- We prove the equalities in the second part of this Theorem.
  This is straightforward to see that $\mathcal{A}_{\mathbb{E}_{C_p,1}} \subseteq \mathcal{A}_{\mathfrak{J}_1}$, $\mathcal{A}_{\mathbb{E}_{C_p^2,C_p^2}} \subseteq \mathcal{A}_{\mathfrak{J}_2}$ and $\mathcal{A}_{\mathbb{E}_{C_p^3,C_p^2}} \subseteq \mathcal{A}_{\mathfrak{J}_3}$.
  Let us show the reverse inclusions. Set $\bar{T} = T/\Phi(T)$ and $\bar{S} = S\Phi(T)/\Phi(T)$.
  1. Let $(T, S)$ be a slice such that $T \neq S$ i.e. $(T, S) \in \mathcal{A}_{\mathfrak{J}_1}$.
     If $S$ is cyclic then there exists a normal subgroup $N$ of the group $T$ of index $p$ such that $S \subseteq N$. Moreover,

     $$T/N \simeq C_p \text{ and } SN/N = 1.$$

     Since the slice $(T/N, SN/N)$ belongs to $\mathcal{A}_{\mathbb{E}_{C_p,1}}$, it follows from Condition **C** that the slice $(T, S) \in \mathcal{A}_{\mathbb{E}_{C_p,1}}$.
     If $S$ is noncyclic then by [Proposition 8.5](), we have $(E \times \bar{T}), E \times \bar{S}) \in \mathcal{A}_{\mathfrak{J}_1}$ whenever $E \times \bar{S}$ is noncyclic. Set $E = C_p \times C_p$, then the slice

     $$(E \times \bar{T}, E \times \bar{S}) = (C_p, 1) \times (C_p \times \bar{T}, E \times \bar{S})$$

     lies in $\mathcal{A}_{\mathbb{E}_{C_p,1}}$. Thus the same argument show that $(T, S) \in \mathcal{A}_{\mathbb{E}_{C_p,1}}$.
  2. Assume that $(T, S)$ is a slice where the group $S$ is noncyclic then for any elementary abelian group such that $|E| \geq p^2$, we have $(T, S) \in \mathcal{A}_{\mathfrak{J}_2}$ if and only if $(E \times \bar{T}, E \times \bar{S}) \in \mathcal{A}_{\mathfrak{J}_2}$.
     In particular, we can set $E = C_p^2$ and [Corollary 7.3]() implies that

     $$(E \times \bar{T}, E \times \bar{S}) \in \mathcal{A}_{\mathbb{E}_{C_p^2,C_p^2}}.$$

     Hence $(T, S) \in \mathcal{A}_{\mathbb{E}_{C_p^2,C_p^2}}$.
  3. Assume that $S$ is noncyclic, $T \neq S$ and $(T, S) \in \mathcal{A}_{\mathfrak{J}_3}$ then $(C_p^2 \times \bar{T}, C_p^2 \times \bar{S}) \in \mathcal{A}_{\mathfrak{J}_3}$. There are integers $n > 1$ and $r < n$ such that

     $$\bar{T} = C_p^n \text{ and } \bar{S} = C_p^r.$$

     Thus

     $$(C_p^2 \times \bar{T}, C_p^2 \times \bar{S}) = (C_p^3, C_p^2) \times (C_p^{n-1}, C_p^r) \in \mathcal{A}_{\mathbb{E}_{C_p^3,C_p^2}},$$

     by [Corollary 7.3](). □

**Remark 8.10.** There is a unital Green biset functor homomorphism $i : \mathbb{K}B_p \to \mathbb{K}\Xi_p$ defined at a $p$-group $G$ by sending the (class of a) finite $G$-set $X$ to the (class of the) identity morphism of $X$.

It was shown moreover in [4] that $\mathbb{K}B_p$ has a unique proper non zero subfunctor, equal to the kernel $K$ of the linearization morphism $\mathbb{K}B_p \to \mathbb{K}R_{\mathbb{Q},p}$. One can check accordingly that $i(\mathbb{K}B_p) \cap \mathfrak{J}_1 = 0$, while $i(\mathbb{K}B_p) \cap \mathfrak{J}_2 = i(\mathbb{K}B_p) \cap \mathfrak{J}_4 = i(K)$.

## 9. The counterexample

Let $G$ be a finite $p$-group.

Consider $\mathfrak{I}_3(G)$, the ideal of $\mathbb{K}\Xi_p(G)$ generated by elements of the form $\xi_{T,S}^G$ where $(T, S)$ is any slice of $G$ such that $S \neq T$ and $S$ is noncyclic. It follows from Theorem 8.9 that the $\mathbb{K}\Xi_p$-submodule $\mathfrak{I}_3$ of $\mathbb{K}\Xi_p$ is simple and for any finite $p$-group $G$, the dimension of $\mathfrak{I}_3(G)$ is equal to the number of conjugacy classes of slices $(T, S)$ of $G$, such that $S \neq T$ and $S$ noncyclic.

Assume that $\mathfrak{I}_3(G) \neq 0$. Then there exists a slice $(T, S)$ of $G$ such that $S \neq T$ and $S$ noncyclic.

In this case, it follows easily that $S$ has at least $p^2$ elements. Since $S$ is a proper subgroup of $T$, it follows that

$(*)$  $|T| \geq p^3$ and $|G| \geq p^3$.

Now we assume that $|G| = p^3$, then $T = G$ and $S \simeq C_p \times C_p$.

If $p \neq 2$, then $G$ is isomorphic to one of the groups $C_p^3$, $C_{p^2} \times C_p$, $M_{p^3} = C_{p^2} \rtimes C_p$ or $X_{p^3} = Syl_p\big(Gl(3, \mathbb{F}_p)\big)$.

The groups $M_{p^3} = C_{p^2} \rtimes C_p$ and $X_{p^3} = Syl_p\big(Gl(3, \mathbb{F}_p)\big)$ are the two extraspecial $p$-groups of order $p^3$ (see [6]).

The lattices of subgroups of $M_{p^3}$ and $X_{p^3}$ for $p = 3$ are respectively:

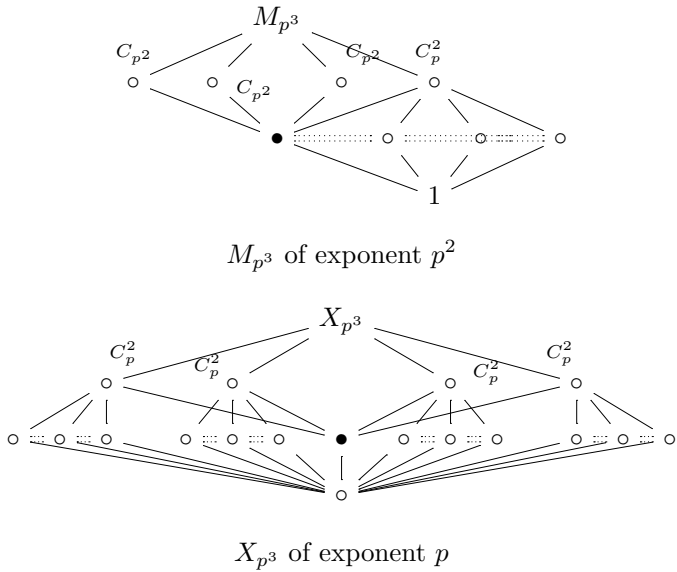

$M_{p^3}$ of exponent $p^2$

$X_{p^3}$ of exponent $p$

An horizontal dotted link between two vertices means that the corresponding subgroups are conjugate. The vertex marked with a ● is the centre of the group.

The lattice of subgroups of the $C_{p^2} \times C_p$ for $p = 3$ is:

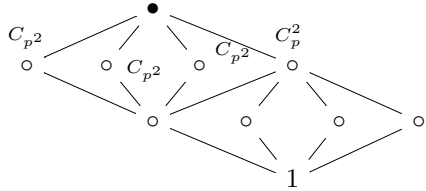

We have therefore established that the slices $(C_p^3, C_p \times C_p)$, $(C_{p^2} \times C_p, C_p \times C_p)$, $(C_{p^2} \rtimes C_p, C_p \times C_p)$ and $(X_{p^3}, C_p \times C_p)$ belong to $\mathcal{A}_{\mathfrak{J}_3}$.

By $(*)$ if $K$ is a finite $p$-group such that $|K| < p^3$ then $\mathfrak{J}_3(K) = 0$.

Thus,

$$\{C_p^3, C_{p^2} \times C_p, M_{p^3}, X_{p^3}\} \subseteq \mathrm{Min}(\mathfrak{J}_3)$$

and

| $G$ | $C_p^3$ | $C_{p^2} \times C_p$ | $M_{p^3}$ | $X_{p^3}$ |
|---|---|---|---|---|
| $\dim \mathfrak{J}_3(G)$ | $p^2 + p + 1$ | 1 | 1 | $p + 1$ |

If $p = 2$, then $G$ is isomorphic to one of the following groups: $C_4 \times C_2$, $C_2^3$, $D_8$.

The lattices of subgroups of $C_4 \times C_2$, $C_2 \times C_2 \times C_2$ and $D_8$ are respectively:

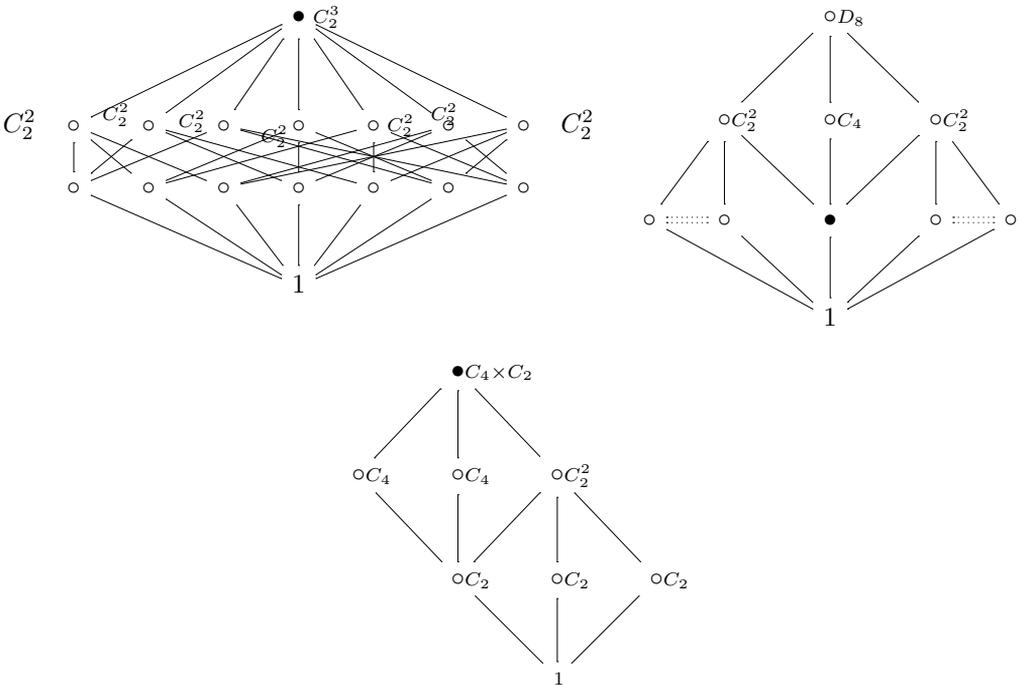

The dimension of $\mathfrak{I}_3(G)$ is therefore the following:

| $G$ | $C_2^3$ | $C_4 \times C_2$ | $D_8$ |
|---|---|---|---|
| dim $\mathfrak{I}_3(G)$ | 7 | 1 | 2 |

## Acknowledgments

The author is grateful to Serge Bouc (university of Picardie-Jules Verne) and Pr. Oumar Diankha (university of Cheikh Anta Diop) for their comments and guidance.

## References

[1] Serge Bouc, Foncteurs d'ensembles munis d'une double action, J. Algebra 183 (0238) (2008) 5067–5087.

[2] Serge Bouc, Bisets Functors for Finite Groups, Lecture Notes in Math., vol. 1990, Springer, 2010.

[3] Serge Bouc, The slice Burnside ring and the section Burnside ring of a finite group, Compos. Math. 148 (1) (2012) 868–906.

[4] Serge Bouc, Jacques Thévenaz, The group of endo-permutation modules, Invent. Math. 139 (2000) 275–349.

[5] Maxime Ducellier, A study of a simple $p$-permutation functor, J. Algebra 447 (2016) 367–382.

[6] Daniel Gorenstein, Finite Groups, Harper and Row, 1968.

[7] Nadia Romero, Tesis, Universidad Nacional Autónoma de México, 2011, http://132.248.9.195/ ptb2011/noviembre/0674682/Index.html.

[8] Nadia Romero, On fibred biset functors with fibres of order prime and four, J. Algebra 387 (2013) 185–194.

[9] Radu Stancu, Jacques Thévenaz, Serge Bouc, Simple biset functors and double Burnside ring, J. Pure Appl. Algebra 217 (2013) 546–566.